\documentclass[a4paper,12pt]{article}
\usepackage{amssymb,fullpage,amsmath,amsthm,graphicx}
\newtheorem{thm}{Theorem}[section]
\newtheorem{defn}{Definition}[section]
\newtheorem{cor}{Corollary}[section]

\newtheorem{lem}{Lemma}[section]
\newtheorem{prop}{Proposition}[section]

\newtheorem{assu}{Assumption}[section]

\numberwithin{equation}{section}

\begin{document}

\title{Scaling limits for simple random walks\\on random ordered graph trees}\author{David Croydon\footnote{Dept of Statistics,
University of Warwick, Coventry, CV4 7AL, UK;
{d.a.croydon@warwick.ac.uk.}}\\
\tiny{UNIVERSITY OF WARWICK}}
\date{15 March 2010}
\maketitle

\begin{abstract}
Consider a family of random ordered graph trees $(T_n)_{n\geq1}$, where $T_n$ has $n$ vertices. It has previously been established that if the associated search-depth processes converge to the normalised Brownian excursion when rescaled appropriately as $n\rightarrow \infty$, then the simple random walks on the graph trees have the Brownian motion on the Brownian continuum random tree as their scaling limit. Here, this result is extended to demonstrate the existence of a diffusion scaling limit whenever the volume measure on the limiting real tree is non-atomic, supported on the leaves of the limiting tree, and satisfies a polynomial lower bound for the volume of balls. Furthermore, as an application of this generalisation, it is established that the simple random walks on a family of Galton-Watson trees with a critical infinite variance offspring distribution, conditioned on the total number of offspring, can be rescaled to converge to the Brownian motion on a related $\alpha$-stable tree.
\end{abstract}

\section{Introduction}

If $(T_n)_{n\geq 1}$ is a family of random ordered graph trees, where $T_n$ has $n$ vertices for each $n$, and the rescaled graph trees $n^{-1/2}T_n$ converge suitably to the Brownian continuum random tree, then the associated simple random walks can be rescaled to converge to a diffusion limit, namely the Brownian motion on the Brownian continuum random tree, see \cite{Croydoncbp}, Theorem 1.2. (Note that, in \cite{Croydoncbp}, the Brownian continuum random tree was referred to simply as the continuum random tree. The `Brownian' is included here to make clear the distinction between this specific random real tree and the other random real trees that will feature in the discussion.) The collection of random graphs to which this result applies includes the case when $T_n$ is a Galton-Watson tree with an offspring distribution, $(p_i)_{i\geq0}$ say, that is aperiodic (in particular, it is not supported on a proper subgroup of $\mathbb{Z}$), critical (mean one), and has finite variance, conditioned on its total progeny being equal to $n$. Motivated by the problem of extending the methods of \cite{Croydoncbp} to deal with the case when the finite variance assumption is dropped, in this article we prove convergence results for the simple random walks on a much broader class of graph trees. Although the overall structure of the argument used here closely matches that of \cite{Croydoncbp}, in that article several rather precise properties of the Brownian continuum random tree were applied, meaning that extensions to more general limiting trees could not be made immediately. Here, new techniques are introduced to allow us to remove some of the restrictions of \cite{Croydoncbp}, and in particular prove a scaling limit result for the simple random walks on critical Galton-Watson trees with infinite variance offspring distributions.

Denote the Brownian continuum random tree and its canonical measure by $(\mathcal{T},\mu)$ (see \cite{Aldous2} for background). Let $\rho$ be a distinguished vertex of $\mathcal{T}$ and call it the root of $\mathcal{T}$. Suppose $\mathcal{T}(k)$ is the minimal subtree of $\mathcal{T}$ spanning $\rho$ and a $k$-sample of $\mu$-random vertices; such sets can be constructed so that $\mathcal{T}(k)\subseteq\mathcal{T}(k+1)$, and also we can assume that $\cup_k\mathcal{T}(k)$ is dense in $\mathcal{T}$, because $\mu$ has full-support, almost-surely. The set $\mathcal{T}(k)$ consists of a finite number of finite length line segments and therefore we can well-define the one-dimensional Hausdorff measure on $\mathcal{T}(k)$ normalised to be a probability measure, $\lambda^{(k)}$ say. In \cite{Croydoncbp}, the fact that $\lambda^{(k)}$ converges weakly to $\mu$ as probability measures on $\mathcal{T}$ was used in a time-change argument to show that certain Brownian motions on the subtrees $\mathcal{T}(k)$ converge to the Brownian motion on the Brownian continuum random tree (\cite{Croydoncbp}, Lemma 3.1). Whilst this assumption is known to hold for the Brownian continuum random tree, it has not been established for the related $\alpha$-stable trees, $\alpha\in(1,2)$, that are the scaling limits of infinite variance Galton-Watson trees (see \cite{rrt}, Section 4, for a brief introduction to $\alpha$-stable trees and this convergence result). By considering alternative Brownian motions on $\mathcal{T}(k)$, we are able to show that the assumption that $\lambda^{(k)}$ converges to $\mu$ is redundant (Proposition \ref{prockconv} is the relevant convergence result). To proceed as we do, however, we are required to check that the time-change additive functionals we consider satisfy a tightness property (Lemma \ref{akprops}), and our proof of this relies on the technical result that the local times of the Brownian motion on the limiting tree diverge uniformly (Lemma \ref{ltdiverge}).

A second fact about the Brownian continuum random tree applied in \cite{Croydoncbp} is that it has no vertex of degree greater than three, by which we mean that the number of connected components of $\mathcal{T}\backslash\{\sigma\}$ is no greater than three for any $\sigma\in\mathcal{T}$. This property was important because it meant that the graph trees $T_n(k)$ constructed analogously to $\mathcal{T}(k)$ eventually had to have the same `tree-shape' as $\mathcal{T}(k)$ as $n$ became large. However, $\alpha$-stable trees with $\alpha\in(1,2)$ admit infinite branch-points (\cite{LegallDuquesne}, Theorem 4.6), and so the same argument does not apply in general. To overcome this problem, in Section \ref{ell1} we establish a method to demonstrate the convergence of Brownian motions on a sequence of finite-branched trees that converge in a specific way to the subtrees $\mathcal{T}(k)$ when the limiting tree $\mathcal{T}$ has branch-points of arbitrary degree.
%\footnote{Let us remark that in \cite{Croydoncbp} there was actually an error in the proof of the convergence result for Brownian motions on finite subtrees that appeared as equation (32). Since this is corrected by results in Section \ref{ell1} of this article, we do not explain the problem or how to correct it in a more elementary fashion.}.
Our argument involves considering approximations to the processes of interest that jump over branch-points, so that the precise geometry of the trees at branch-points is not seen by the approximations.

Specifically, the limit space $\mathbb{T}^*$ that we study in this article is the collection of pairs $(\mathcal{T},\mu)$, where $\mathcal{T}=(\mathcal{T},d_\mathcal{T})$ is a compact real tree (see \cite{LegallDuquesne}, Definition 2.1, for example) and $\mu$ is a non-atomic Borel probability measure on $\mathcal{T}$ that satisfies
\begin{equation}\label{mulower}
\liminf_{r\rightarrow 0}\frac{\inf_{\sigma\in\mathcal{T}}\mu(B(\sigma,r))}{r^\kappa}>0,
\end{equation}
for some $\kappa>0$, where $B(\sigma,r)$ is the open ball of radius $r$ (with respect to the metric $d_\mathcal{T}$) centred at $\sigma\in\mathcal{T}$. Furthermore, for $(\mathcal{T},\mu)\in\mathbb{T}^*$, we assume that $\mu$ is supported on the leaves of $\mathcal{T}$, where by saying that $\sigma\in\mathcal{T}$ is a leaf of $\mathcal{T}$ we mean that $\mathcal{T}\backslash\{\sigma\}$ is connected. We note that the Brownian continuum random tree is an element of $\mathbb{T}^*$ almost-surely, with (\ref{mulower}) being satisfied whenever $\kappa>2$, see \cite{Croydoncrt}, Theorem 1.2. Furthermore, we will later check that all the above properties are also almost-surely satisfied by $\alpha$-stable trees for $\alpha\in(1,2)$, with (\ref{mulower}) holding for any choice of $\kappa>\alpha/(\alpha-1)$.

The main result of this article is the following theorem, which is stated within the framework of \cite{Croydoncbp}, so that: $\ell^1$ is the Banach space of infinite sequences of real numbers equipped with the metric $d_{\ell^1}$ induced by the norm $\sum_{i\geq 1}|x(i)|$ for $x\in\ell^1$; $\mathcal{K}(\ell^1)$ is the space of compact subsets of $\ell^1$ equipped with the Hausdorff topology; $\mathcal{M}_1(\ell^1)$ is the space of Borel probability measures on $\ell^1$ equipped with the topology of weak convergence; and $\mathcal{M}_1(C([0,1],\ell^{1}))$ is the space of Borel probability measures on $C([0,1],\ell^{1})$, also equipped with the topology of weak convergence. The set $\mathcal{W}$ is the collection of continuous functions $w:[0,1]\rightarrow \mathbb{R}_+$ such that $w(t)>0$ if and only if $t\in (0,1)$. For an excursion $w\in\mathcal{W}$: $\mathcal{T}_w$ is the rooted real tree associated with $w$; $\mu_w$ is the natural measure on $\mathcal{T}_w$; and $\mathbf{P}^{\mathcal{T}_w,\mu_w}_\rho$ is the law of the Brownian motion on $(\mathcal{T}_w,\mu_w)$ started from the root $\rho=\rho(\mathcal{T}_w)$ (see the end of Section \ref{bmrt} for details). For an ordered graph tree $T_n$ with $n$ vertices, we write $\mu_n$ to represent the uniform probability measure on the vertices of $T_n$, and denote by $\mathbf{P}^{{T}_n}_\rho$ the law of the discrete time simple random walk on ${T}_n$, started from the root, $\rho=\rho({T}_n)$, of ${T}_n$. The search-depth process $w_n$ of $T_n$ is defined in Section \ref{discconv}.

{\thm \label{qthm} Let $(\alpha_n)_{n\geq 1}$ be a positive divergent sequence such that $\alpha_n=o(n)$. If $({T}_n)_{n\geq 1}$ is a sequence of ordered graph trees whose search-depth functions $(w_n)_{n\geq 1}$ satisfy
\[\alpha_n^{-1}w_n\rightarrow w\]
in $C([0,1],\mathbb{R}_+)$ for some $w\in\mathcal{W}$ with $(\mathcal{T}_w,\mu_w)\in\mathbb{T}^*$, then there exists, for each $n$, an isometric embedding $(\tilde{{T}}_n,\tilde{\mu}_n,\tilde{\mathbf{P}}^{{T}_n}_\rho)$ of the triple $({{T}}_n,{\mu}_n,{\mathbf{P}}^{{T}_n}_\rho)$ into $\ell^1$ such that
\[\left(\alpha_n^{-1}\tilde{T}_n,\tilde{\mu}_n(\alpha_n\cdot),{\mathbf{P}}^{{T}_n}_\rho(\{f\in C([0,1],\ell^1):\:(\alpha_n^{-1}f(tn\alpha_n))_{t\in[0,1]}\in\cdot\})\right)\]
converges to $(\tilde{\mathcal{T}},\tilde{\mu},\tilde{\mathbf{P}}^{\mathcal{T},\mu}_\rho)$ in the space $\mathcal{K}(\ell^1)\times\mathcal{M}_1(\ell^1)\times \mathcal{M}_1(C([0,1],\ell^1))$, where the law $\tilde{\mathbf{P}}^{{T}_n}_\rho$ is extended to an element of $\mathcal{M}_1(C([0,1],\ell^1))$ by linear interpolation of discrete time processes, and  $(\tilde{\mathcal{T}},\tilde{\mu},\tilde{\mathbf{P}}^{\mathcal{T},\mu}_\rho)$ is an isometric embedding of $(\mathcal{T}_w,{\mu}_w,{\mathbf{P}}^{\mathcal{T}_w,\mu_w}_\rho)$ into $\ell^1$.}
\bigskip

As an extension to \cite{Croydoncbp}, Theorem 1.2, a random version of this result can be formulated by applying the fact that the isometrically embedded triple $(\tilde{\mathcal{T}},\tilde{\mu},\tilde{\mathbf{P}}^{\mathcal{T},\mu}_\rho)$ can be constructed as a measurable function of a pair $(w,u)$, where $w$ is the relevant excursion and $u$ is an element of $[0,1]^{\mathbb{N}}$. In particular, suppose that $(W,U)$ is a $\mathcal{W}\times[0,1]^\mathbb{N}$-valued random variable, built on a complete probability space with probability measure $\mathbf{P}$, such that $W$ is a random excursion satisfying $\mathbf{P}((\mathcal{T}_W,\mu_W)\in\mathbb{T}^*)=1$, and $U=(U_i)_{i\geq 1}$ is an independent sequence of independent $U[0,1]$ random variables. We can define a probability law $\mathbb{P}$ on $\mathcal{K}(l^1)\times\mathcal{M}_1(\ell^1)\times C([0,1],\ell^1)$ that satisfies
\begin{equation}\label{pdef}
\mathbb{P}\left(A\times B \times C\right)=\int_{C([0,1],\mathbb{R}_+)\times [0,1]^\mathbb{N}} \mathbf{P}((W,U)\in (dw,du))\: \mathbf{1}_{\{\tilde{\mathcal{T}}\in A,\:\tilde{\mu}\in B\}}\tilde{\mathbf{P}}^{\mathcal{T},\mu}_\rho(C),
\end{equation}
for every measurable $A\subseteq\mathcal{K}(\ell^1)$, $B\subseteq\mathcal{M}_1(\ell^1)$, and $C\subseteq C([0,1],\ell^1)$. In the discrete case, as in \cite{Croydoncbp}, let $({T}_n)_{n\geq 1}$ be a sequence of random ordered graph trees with corresponding search-depth functions $(W_n)_{n\geq 1}$, and suppose these are built on our underlying probability space independently of the random variable $U$. We can assume that $(\tilde{{T}}_n,\tilde{\mu}_n,\tilde{\mathbf{P}}^{{T}_n}_\rho)$ are constructed measurably from $(W_n,U)$. Moreover, there exists a unique probability law $\mathbb{P}_n$ on $\mathcal{K}(\ell^1)\times\mathcal{M}_1(\ell^1)\times C([0,1],\ell^1)$ that satisfies
\begin{equation}\label{pndef}
\mathbb{P}_n\left(A\times B \times C\right)=\int_{C([0,1],\mathbb{R}_+)\times [0,1]^\mathbb{N}} \mathbf{P}((W_n,U)\in (dw,du))\: \mathbf{1}_{\{\tilde{{T}}_n\in A,\:\tilde{\mu}_n\in B\}}\tilde{\mathbf{P}}^{{T}_n}_\rho(C),
\end{equation}
for every measurable $A\subseteq\mathcal{K}(\ell^1)$, $B\subseteq\mathcal{M}_1(\ell^1)$, and $C\subseteq C([0,1],\ell^1)$. In the following result, the rescaling operator $\Theta_n$ is defined on $\mathcal{K}(\ell^1)\times \mathcal{M}_1(\ell^1)\times C([0,1],\ell^{1})$ so that if $(\tilde{K},\tilde{\nu},\tilde{f})$ is an element of this space, then
\[\Theta_n(\tilde{K},\tilde{\nu},\tilde{f}):=(\alpha_n^{-1}\tilde{K},\tilde{\nu}(\alpha_n\cdot),(\alpha_n^{-1}\tilde{f}(tn\alpha_n))_{t\in[0,1]}).\]

{\thm\label{athm}  Let $(\alpha_n)_{n\geq 1}$ be a positive divergent sequence such that $\alpha_n=o(n)$. Suppose that $({T}_n)_{n\geq 1}$ is a sequence of random ordered graph trees whose rescaled search-depth functions $(\alpha_n^{-1}W_n)_{n\geq 1}$ converge in distribution to $W$ in $C([0,1],\mathbb{R}_+)$, where $W$ is a random excursion satisfying $\mathbf{P}((\mathcal{T}_W,\mu_W)\in\mathbb{T}^*)=1$. If $\mathbb{P}$ and $(\mathbb{P}_n)_{n\geq 1}$ are the unique probability measures satisfying (\ref{pdef}) and (\ref{pndef}), respectively, then
\[\mathbb{P}_n\circ\Theta_n^{-1}\rightarrow \mathbb{P}\]
weakly as measures on the space $\mathcal{K}(\ell^1)\times\mathcal{M}_1(\ell^1)\times C([0,1],\ell^1)$.}
\bigskip

In addition to the above convergence results for the laws of the simple random walks on the graph trees $(T_n)_{n\geq 1}$, let us remark that, in the settings of Theorem \ref{qthm} or \ref{athm}, it is possible to proceed as in \cite{CHLLT}, Section 7.2, to deduce related local limit theorems demonstrating that the associated discrete transition densities can be rescaled to converge in an appropriate space to the transition densities of the Brownian motion on the limiting space.

To prove the results stated above, we commence, in Section \ref{bmrt}, by proving some general results about Brownian motion and the corresponding local times on real trees. In Section \ref{bmrt}, we also describe the procedures we use for embedding real trees into $\ell^1$, and the connection between real trees and excursions. The heart of the article is Section \ref{ell1}, which is where we establish some key results about the convergence of Brownian motions on real trees with a finite number of branches embedded into $\ell^1$. Section \ref{discconv}, which contains the proofs of Theorems \ref{qthm} and \ref{athm}, explains how the argument of \cite{Croydoncbp} can be reworked to apply in our more general setting. Finally, in Section \ref{stable}, we describe the application of our results to Galton-Watson trees with infinite variance offspring distributions and the related $\alpha$-stable trees.

\section{Brownian motion and local times on real trees}\label{bmrt}

Suppose $\mathcal{T}=(\mathcal{T},d_\mathcal{T})$ is a compact real tree and $\mu$ is a finite Borel measure on $\mathcal{T}$ with full support. To avoid trivialities, assume that $\mathcal{T}$ contains more than one point. It is possible (see \cite{Croydoncbp}, Proposition 2.2) to construct a strong Markov process
\[X^{\mathcal{T},\mu}=\left(\left(X^{\mathcal{T},\mu}_t\right)_{t\geq 0}, \mathbf{P}_\sigma^{\mathcal{T},\mu}, \sigma\in \mathcal{T}\right),\]
with continuous sample paths that is reversible with respect to its invariant measure $\mu$ and satisfies the following properties.
\newcounter{listcount}
\begin{list}{\roman{listcount})}{\usecounter{listcount} \setlength{\rightmargin}{\leftmargin}}
\item For $\sigma_1,\sigma_2\in \mathcal{T}$, $\sigma_1\neq \sigma_2$, we have
\[\mathbf{P}_{\sigma}^{\mathcal{T},\mu}\left(h({\sigma_1})<h({\sigma_2}) \right)=\frac{d_\mathcal{T}(b^\mathcal{T}(\sigma,\sigma_1,\sigma_2),\sigma_2)}{d_\mathcal{T}(\sigma_1,\sigma_2)},\hspace{20pt}\forall \sigma\in \mathcal{T},\]
where $h({\sigma}):=\inf\{t>0:\:X^{\mathcal{T},\mu}_t=\sigma\}$ is the hitting time of $\sigma\in \mathcal{T}$, and $b^\mathcal{T}(\sigma,\sigma_1,\sigma_2)$ is the unique branch-point of $\sigma$, $\sigma_1$ and $\sigma_2$ in $\mathcal{T}$. In particular, if $[[\sigma, \sigma_1]]$, $[[\sigma_1, \sigma_2]]$ and $[[\sigma_2, \sigma]]$ are the unique injective paths between the relevant pairs of vertices, then $b^\mathcal{T}(\sigma,\sigma_1,\sigma_2)$ is
the unique point in the set $[[\sigma, \sigma_1]]\cap[[\sigma_1, \sigma_2]]\cap[[\sigma_2, \sigma]]$.

\item For $\sigma_1,\sigma_2\in \mathcal{T}$, the mean occupation measure for the process started at $\sigma_1$ and killed on hitting $\sigma_2$ has density $2d_\mathcal{T}(b^\mathcal{T}(\sigma,\sigma_1,\sigma_2),\sigma_2)\mu(d\sigma)$, so that
    \[\mathbf{E}_{\sigma_1}^{\mathcal{T},\mu}\int_0^{h(\sigma_2)}f(X_s)ds=2\int_{\mathcal{T}}f(\sigma) d_\mathcal{T}(b^\mathcal{T}(\sigma,\sigma_1,\sigma_2),\sigma_2)\mu(d\sigma),\]
    for every continuous bounded function $f:\mathcal{T}\rightarrow \mathbb{R}$.
\end{list}
In the terminology of \cite{Aldous2}, Section 5.2, $X^{\mathcal{T},\mu}$ is Brownian motion on $(\mathcal{T},\mu)$, and is in fact uniquely determined by these properties. Moreover, that $X^{\mathcal{T},\mu}$ admits jointly measurable local times $(L_t(\sigma))_{\sigma\in\mathcal{T},t\geq0}$ can be checked as in the proof of \cite{Croydoncrt}, Lemma 8.2. In the arguments of subsequent sections we will require further that $(L_t(\sigma))_{\sigma\in\mathcal{T},t\geq0}$ is jointly continuous in $t$ and $\sigma$ and we will demonstrate that this is the case whenever $\mu$ satisfies a polynomial lower bound of the form of (\ref{mulower}). In the proof of this result, we apply the two following properties that (\ref{mulower}) implies.

{\lem \label{est} Suppose $\mathcal{T}$ is a compact real tree and $\mu$ is a finite Borel measure on $\mathcal{T}$ that satisfies (\ref{mulower}) for some $\kappa>0$.\\
(a) If $N(\mathcal{T},r)$ is the smallest number of balls of radius $r$ needed to cover $\mathcal{T}$, then
\[\limsup_{r\rightarrow 0}r^\kappa{N(\mathcal{T},r)}<\infty.\]
(b) The Markov process $X^{\mathcal{T},\mu}$ admits a transition density $(p_t(\sigma,\sigma'))_{\sigma,\sigma'\in\mathcal{T},t>0}$ that satisfies
\[\limsup_{t\rightarrow 0}{t^{\frac{\kappa}{\kappa+1}}}{\sup_{\sigma,\sigma'\in\mathcal{T}}p_t(\sigma,\sigma')}<\infty.\]}
\begin{proof} The proof of (a) is elementary. Part (b) can be obtained by applying a general heat kernel bound of the type proved in \cite{Kumagai}, Proposition 4.1, or \cite{Croydon}, Proposition 5, for example.
\end{proof}

{\lem \label{ltcont} If $\mathcal{T}$ is a compact real tree and $\mu$ is a finite Borel measure on $\mathcal{T}$ that satisfies (\ref{mulower}) for some $\kappa>0$, then the local times $(L_t(\sigma))_{\sigma\in\mathcal{T},t\geq0}$ of $X^{\mathcal{T},\mu}$ are jointly continuous in $t$ and $\sigma$, $\mathbf{P}^{\mathcal{T},\mu}_{\sigma'}$-a.s., for every $\sigma'\in\mathcal{T}$.}
\begin{proof} Given the estimates of Lemma \ref{est}, the proof is identical to that of \cite{Croydoncbp}, Lemma 2.5. In particular, it is easily checked from Lemma \ref{est}(b) that the 1-potential density $u(\sigma,\sigma'):=\int_{0}^\infty e^{-t}p_t(\sigma,\sigma')$ is finite for all $\sigma,\sigma'\in\mathcal{T}$. As a consequence, by applying \cite{MarcusRosen}, Theorem 1, we see that the continuity of local times is equivalent to the continuity of the centred Gaussian process $(G(\sigma))_{\sigma\in \mathcal{T}}$ with covariances given by $(u(\sigma,\sigma'))_{\sigma,\sigma'\in\mathcal{T}}$. By \cite{DudleyGauss}, Theorem 2.1, for the latter process to be continuous, it is enough that the integral $\int_0^1 \sqrt{\ln N(\mathcal{T},r)}dr$ is finite, which in view of Lemma \ref{est}(a) is clearly the case.
\end{proof}

The local times of $X^{\mathcal{T},\mu}$ will be used in a time-change argument that depends on their uniform divergence, which we prove now. Note that, by definition, the Brownian motion on $(\mathcal{T},\mu)$ satisfies $\mathbf{E}^{\mathcal{T},\mu}_\sigma h(\sigma')\leq2\mathrm{diam}(\mathcal{T})\mu(\mathcal{T})<\infty$ for every $\sigma,\sigma'\in\mathcal{T}$, where $\mathrm{diam}(\mathcal{T})$ is the diameter of the metric space $(\mathcal{T},d_\mathcal{T})$.

{\lem \label{ltdiverge} Suppose $\mathcal{T}$ is a compact real tree and $\mu$ is a finite Borel measure on $\mathcal{T}$ that satisfies (\ref{mulower}) for some $\kappa>0$. For every $\sigma\in\mathcal{T}$, $\mathbf{P}_\sigma^{\mathcal{T},\mu}$-a.s. we have
\[\lim_{t\rightarrow\infty}\inf_{\sigma'\in\mathcal{T}}L_t(\sigma')=\infty.\]}
\begin{proof} Fix $\sigma,\sigma'\in\mathcal{T}$. By \cite{MarcusRosen}, Lemma 3.6, we have that $\mathbf{E}^{\mathcal{T},\mu}_{\sigma'}\int_0^\infty e^{-t}d_tL_t(\sigma')>0$. Hence there is a strictly positive $\mathbf{P}^{\mathcal{T},\mu}_{\sigma'}$-probability that $L_t(\sigma')>0$ for large $t$. Applying the joint continuity of the local times, it follows that there exist $r=r(\sigma')>0$, $\varepsilon=\varepsilon(\sigma')>0$ and $t_0=t_0(\sigma')<\infty$ such that
\begin{equation}\label{posprob}
\mathbf{P}^{\mathcal{T},\mu}_{\sigma'}\left(\inf_{\sigma''\in B(\sigma',r)}L_{t_0}(\sigma'')> \varepsilon\right)>0.
\end{equation}
Now, set $h(\sigma',t_0,\sigma):=\inf\{t>t_0+h(\sigma'):X^{\mathcal{T},\mu}_t=\sigma\}$.
Applying the observation made above this lemma about the finite moments of hitting times, the strong Markov property and (\ref{posprob}), it is easy to check that $h(\sigma',t_0,\sigma)$ is finite, $\mathbf{P}^{\mathcal{T},\mu}_\sigma$-a.s., and also
\begin{equation}\label{posprob2}
\mathbf{P}^{\mathcal{T},\mu}_\sigma\left(\inf_{\sigma''\in B(\sigma',r)}L_{h(\sigma',t_0,\sigma)}(\sigma'')> \varepsilon\right)>0.
\end{equation}
Observe that the additivity of local times and the strong Markov property implies that
\[\liminf_{t\rightarrow\infty}\inf_{\sigma''\in B(\sigma',r)}L_{t}(\sigma'')\geq \sum_{i=1}^\infty \xi_i,\]
where, under $\mathbf{P}_\sigma^{\mathcal{T},\mu}$, $(\xi_i)_{i=1}^\infty$ are independent copies of $\inf_{\sigma''\in B(\sigma',r)}L_{h(\sigma',t_0,\sigma)}(\sigma'')$. The strong law of large numbers lets it be deduced from (\ref{posprob2}) that the right-hand side of the above inequality is infinite, $\mathbf{P}_\sigma^{\mathcal{T},\mu}$-a.s., which proves the uniform divergence of local times uniformly over $B(\sigma',r)$, $\mathbf{P}_\sigma^{\mathcal{T},\mu}$-a.s.

To extend the conclusion of the previous paragraph, note that $(B(\sigma',r(\sigma'))_{\sigma'\in\mathcal{T}}$ is an open cover for $\mathcal{T}$. Thus, by the compactness of $\mathcal{T}$, it admits a finite subcover, $(B(\sigma_i,r(\sigma_i))_{i=1}^N$, and clearly
\[\lim_{t\rightarrow\infty}\inf_{\sigma'\in\mathcal{T}}L_t(\sigma')= \min_{i=1,\dots,N}\lim_{t\rightarrow\infty}\inf_{\sigma'\in B(\sigma_i,r(\sigma_i))}L_{t}(\sigma').\]
Since, by our above argument, the right-hand side of this expression is infinite, $\mathbf{P}_\sigma^{\mathcal{T},\mu}$-a.s., the proof is complete.
\end{proof}

Following \cite{Croydoncbp}, we now show how $X^{\mathcal{T},\mu}$ can be approximated by a family of Brownian motions on subtrees of $\mathcal{T}$ with a finite number of branches. Henceforth, we suppose that the real tree $\mathcal{T}$ has a distinguished vertex $\rho\in\mathcal{T}$ called the root and consider a dense sequence of vertices $(\sigma_i)_{i=1}^\infty$ in $\mathcal{T}$. Without loss of generality, we assume that $(\sigma_i)_{i=1}^\infty$ are distinct and $\sigma_i\neq \rho$ for any $i$. For each $k\geq 1$, define a subset of $\mathcal{T}$ by
\begin{equation}\label{Tkdef}
\mathcal{T}(k):=\bigcup_{i=1}^k[[\rho,\sigma_i]],
\end{equation}
where, as above,  for $\sigma,\sigma'\in\mathcal{T}$, $[[\sigma,\sigma']]$ is the unique injective path in $\mathcal{T}$ from $\sigma$ to $\sigma'$. Clearly $\mathcal{T}(k)$ is a compact real tree when endowed with the appropriate restriction of $d_\mathcal{T}$, and we set its root to be $\rho$, which is contained in $\mathcal{T}(k)$ by construction. The natural projection $\phi_{\mathcal{T},\mathcal{T}(k)}$ from $\mathcal{T}$ to $\mathcal{T}(k)$ is obtained by setting $\phi_{\mathcal{T},\mathcal{T}(k)}(\sigma)$ to be the unique point in $\mathcal{T}(k)$ satisfying
\begin{equation}\label{projproj}
d_\mathcal{T}(\sigma,\phi_{\mathcal{T},\mathcal{T}(k)}(\sigma))=\inf_{\sigma'\in\mathcal{T}(k)}d_\mathcal{T}(\sigma,\sigma').
\end{equation}
It is elementary to check that $\phi_{\mathcal{T},\mathcal{T}(k)}$ is continuous and $\sup_{x\in\mathcal{T}}d_\mathcal{T} (x,\phi_{\mathcal{T},\mathcal{T}(k)}(x))\rightarrow 0$ (cf. \cite{Croydoncbp}, Lemma 2.4). Consequently, $\mu^{(k)}:=\mu\circ\phi_{\mathcal{T},\mathcal{T}(k)}^{-1}$ defines a Borel probability measure on $\mathcal{T}(k)$ with full support for each $k$, and $\mu^{(k)}\rightarrow \mu$ weakly as probability measures on $\mathcal{T}$. Since $\mathcal{T}(k)$ is a compact real tree containing more than one point and $\mu^{(k)}$ is a Borel probability measure on $\mathcal{T}(k)$ with full support, we can construct the Brownian motion $X^{\mathcal{T}(k),\mu^{(k)}}$ on $(\mathcal{T}(k),\mu^{(k)})$. For these processes, we are able to deduce the following convergence result. Since it is a relatively simple adaptation of \cite{Croydoncbp}, Lemma 3.1, we only sketch the proof.

{\prop \label{prockconv} If $(\mathcal{T},\mu)$ and $\{(\mathcal{T}(k),\mu^{(k)})\}_{k=1}^\infty$ are as described above, then
\[\mathbf{P}_{\rho}^{\mathcal{T}(k),\mu^{(k)}}\rightarrow \mathbf{P}_\rho^{\mathcal{T},\mu}\]
weakly as probability measures on $C(\mathbb{R}_+,\mathcal{T})$.}
\begin{proof} Applying the weak convergence of $\mu^{(k)}$ to $\mu$ and the joint continuity of the local times of $X^{\mathcal{T},\mu}$ (see Lemma \ref{ltcont}), we obtain for every $t\geq 0$ that, $\mathbf{P}_\rho^{\mathcal{T},\mu}$-a.s.,
\[\tilde{A}^{(k)}_t:=\int_{\mathcal{T}(k)}L_t(\sigma)\mu^{(k)}(d\sigma)\rightarrow t.\]
Moreover, an elementary monotonocity argument yields this convergence result uniformly on compact intervals. As a consequence of this, $\tilde{\tau}^{(k)}(t):=\inf\{s:\tilde{A}^{(k)}_s>t\}\rightarrow t$ uniformly on compacts, $\mathbf{P}_\rho^{\mathcal{T},\mu}$-a.s. Now, the trace theorem for Dirichlet forms (see \cite{FOT}, Theorem 6.2.1, for example) allows one to check that the law of $(X_{\tilde{\tau}^{(k)}(t)}^{\mathcal{T},\mu})_{t\geq0}$ under $\mathbf{P}_\rho^{\mathcal{T},\mu}$ is precisely $\mathbf{P}_{\rho}^{\mathcal{T}(k),\mu^{(k)}}$ (cf. \cite{Croydoncbp}, Lemma 2.6), and hence the result follows.
\end{proof}

We continue by presenting a characterisation of $X^{\mathcal{T}(k),\mu^{(k)}}$ as a time-change of another Brownian motion on $\mathcal{T}(k)$. For $k\geq 1$, let $\lambda^{(k)}$ be the one-dimensional Hausdorff measure on $\mathcal{T}(k)$ normalised to have total mass equal to one. Since $\mathcal{T}(k)$ consists of a finite number of line segments, $\lambda^{(k)}$ is a Borel probability measure on $\mathcal{T}(k)$ with full support. Consequently, the Brownian motion $X^{\mathcal{T}(k),\lambda^{(k)}}$ on $(\mathcal{T}(k),\lambda^{(k)})$ exists. Furthermore, it is elementary to check that $\lambda^{(k)}$ satisfies (\ref{mulower}) with $\kappa=1$ and therefore we can apply Lemma \ref{ltcont} to deduce that $X^{\mathcal{T}(k),\lambda^{(k)}}$ admits jointly continuous local times $(L_t^{(k)}(\sigma))_{\sigma\in\mathcal{T}(k),t\geq 0}$. As in \cite{Croydoncbp}, define a continuous additive functional $\hat{A}^{(k)}=(\hat{A}^{(k)}_t)_{t\geq0}$ by setting
\begin{equation}\label{hatak}
\hat{A}^{(k)}_t:=\int_{\mathcal{T}(k)}L^{(k)}_t(\sigma)\mu^{(k)}(d\sigma),
\end{equation}
and its inverse by
\begin{equation}\label{hattau}
\hat{\tau}^{(k)}(t):=\inf\{s:\hat{A}^{(k)}_s>t\}.
\end{equation}
As with the time-change employed in the proof of the previous result, the following lemma is a relatively straightforward consequence of the trace theorem for Dirichlet forms (see \cite{FOT}, Theorem 6.2.1, for the trace theorem and \cite{Croydoncbp}, Lemma 2.6, for a similar application), and so will be stated without proof.

{\lem \label{timechange} Suppose $(\mathcal{T},\mu)$ and $\{(\mathcal{T}(k),\mu^{(k)})\}_{k=1}^\infty$ are as described above. If the process $X^{\mathcal{T}(k),\lambda^{(k)}}$ has law  $\mathbf{P}^{\mathcal{T}(k),\lambda^{(k)}}_\rho$, then the process
\[\left(X^{\mathcal{T}(k),\lambda^{(k)}}_{\hat{\tau}^{(k)}(t)}\right)_{t\geq 0}\]
has law $\mathbf{P}_\rho^{\mathcal{T}(k),\mu^{(k)}}$.}
\bigskip

We now deduce some simple path properties of $\hat{A}^{(k)}$ that will be useful to us later.

{\lem \label{akprops}If $(\mathcal{T},\mu)$ and $\{(\mathcal{T}(k),\mu^{(k)})\}_{k=1}^\infty$ are as described above, then for every $k\geq 1$, $\mathbf{P}_\rho^{\mathcal{T}(k),\lambda^{(k)}}$-a.s., the functions $\hat{A}^{(k)}$ are continuous and strictly increasing. Moreover, for every $t_0\geq0$,
\begin{equation}\label{diverge}
\lim_{t\rightarrow\infty}\limsup_{k\rightarrow\infty}\mathbf{P}_\rho^{\mathcal{T}(k),\lambda^{(k)}}\left(\hat{A}^{(k)}_{t_0}>t\right)=0.
\end{equation}}
\begin{proof} The continuity of $\hat{A}^{(k)}$ follows from the continuity of $(L^{(k)}_t(\sigma))_{\sigma\in\mathcal{T}(k),t\geq0}$, which was noted above Lemma \ref{timechange}. Hence, it remains to show that $\hat{A}^{(k)}_t$ strictly increases in $t$ and satisfies (\ref{diverge}). We start by showing that $(\hat{A}^{(k)}_t)_{t\geq 0}$ is strictly increasing. The following argument holds $\mathbf{P}^{\mathcal{T}(k),\lambda^{(k)}}_\rho$-a.s. Let $s<t$, then clearly \[\int_{\mathcal{T}(k)}(L^{(k)}_t(x)-L^{(k)}_s(x))\lambda^{(k)}(dx)=t-s>0.\]
Hence there exists an $\varepsilon>0$ and a non-empty open set $A\subseteq\mathcal{T}(k)$ such that $L_t^{(k)}(x)-L^{(k)}_s(x)>\varepsilon$ for $x\in A$. Since $\mu^{(k)}$ has full support, it charges every non-empty open set and so we must therefore have that $\hat{A}^{(k)}_s\leq \hat{A}^{(k)}_t-\varepsilon\mu^{(k)}(A)<\hat{A}^{(k)}_t$, which proves the desired result.

We now prove the tightness result of (\ref{diverge}). Use the local times of $X^{\mathcal{T},\mu}$ to define an additive functional ${A}^{(k)}=({A}^{(k)}_t)_{t\geq0}$ by setting
\[{A}^{(k)}_t:=\int_{\mathcal{T}}L_t(\sigma)\lambda^{(k)}(d\sigma),\]
and its inverse ${\tau}^{(k)}$ similarly to the definition of $\hat\tau^{(k)}$ at (\ref{hattau}). By again applying the trace theorem for Dirichlet forms, it can be deduced that the law of the process $(X^{\mathcal{T},\mu}_{{\tau}^{(k)}(t)})_{t\geq 0}$ under $\mathbf{P}_\rho^{\mathcal{T},\mu}$ is the same as that of the process  $X^{\mathcal{T}(k),\lambda^{(k)}}$ under $\mathbf{P}_\rho^{\mathcal{T}(k),\lambda^{(k)}}$. Similarly to \cite{Croydoncbp}, Lemma 3.4, it follows that, under $\mathbf{P}_\rho^{\mathcal{T},\mu}$, the two-parameter process $(L_{{\tau}^{(k)}(t)}(\sigma))_{\sigma\in\mathcal{T}(k),t\geq 0}$
has the same distribution as $(L^{(k)}_t(\sigma))_{\sigma\in\mathcal{T}(k),t\geq 0}$ under $\mathbf{P}_\rho^{\mathcal{T}(k),\lambda^{(k)}}$. Consequently, to complete the proof it will suffice to demonstrate that, for every $t_0\geq 0$, $\mathbf{P}_\rho^{\mathcal{T},\mu}$-a.s.,
\[\limsup_{k\rightarrow\infty}\int_\mathcal{T}L_{{\tau}^{(k)}(t_0)}(\sigma)\mu^{(k)}(d\sigma)<\infty.\]
Recalling that $\mu^{(k)}$ converges weakly to $\mu$, we have that the left-hand side of the above expression is bounded above by
\[\int_\mathcal{T}L_{\sup_{k}{\tau}^{(k)}(t_0)}(\sigma)\mu(d\sigma)=\sup_{k}{\tau}^{(k)}(t_0),\]
and this supremum is finite whenever $\inf_{k}A^{(k)}_t$ diverges as $t\rightarrow\infty$. Applying the uniform divergence of local times proved in Lemma \ref{ltdiverge} and the fact that $\lambda^{(k)}$ is by definition a probability measure, this result holds $\mathbf{P}_\rho^{\mathcal{T},\mu}$-a.s. as required.
\end{proof}

As in \cite{Croydoncbp}, to embed $\mathcal{T}$ into $\ell^1$ we use the sequential embedding of \cite{Aldous3}, Section 2.2. In particular, given a sequence $(\mathcal{T}(k))_{k\geq 1}$ as above it is possible to construct a distance-preserving map $\psi:(\mathcal{T},d_\mathcal{T})\rightarrow (\ell^1, d_{\ell^1})$ that satisfies $\psi(\rho)=0$ and
\begin{equation}\label{proj}
\pi_k(\psi(\sigma))=\psi(\phi_{\mathcal{T},\mathcal{T}(k)}(\sigma))
\end{equation}
for every $\sigma\in\mathcal{T}$ and $k\geq 1$, where $\pi_k$ is the projection map on $\ell^1$ defined by setting $\pi_k(x(1),x(2),\dots)=(x(1),\dots,x(k),0,0,\dots)$. Such a map is determined uniquely by insisting that $\psi(\mathcal{T})\subseteq\{(x(1),x(2),\dots)\in\ell^1:x(i)\geq 0, i=1,2,\dots\}$. We will denote the $\ell^1$-embedded versions of the objects $\mathcal{T},\mu,\mathbf{P}_\rho^{\mathcal{T},\mu},\dots$ by $\tilde{\mathcal{T}},\tilde{\mu},\tilde{\mathbf{P}}_\rho^{\mathcal{T},\mu},\dots$ respectively.

To complete this section, we present the well-known relation between continuous excursions and real trees, and define a collection of pairs of excursions and sequences that will be of interest later in this article. Let $\mathcal{W}$ be defined as in the introduction. For $w\in\mathcal{W}$, define a distance on $[0,1]$ by setting
\begin{equation}\label{dwdef}
d_w(s,t):=w(s)+w(t)-2\inf\{w(r):\:r\in[s\wedge t,s\vee t]\},
\end{equation}
and then use the equivalence $s\sim_w t$ if and only if $d_w(s,t)=0$, to define $\mathcal{T}_w:=[0,1]/\sim_w$. Denoting the canonical projection (with respect to $\sim_w$) from $[0,1]$ to $\mathcal{T}_w$ by $\hat{w}$, it is possible to check that
$d_{\mathcal{T}_w}(\hat{w}(s),\hat{w}(t)):=d_w(s,t)$ defines a metric on
$\mathcal{T}_w$, and also that with this metric $\mathcal{T}_w$ is a compact real tree (see \cite{LegallDuquesne}, Theorem 2.1). The root of the tree $\mathcal{T}_w$ is defined to be the equivalence class $\hat{w}(0)$, and is denoted by $\rho_w$. A Borel probability measure on $\mathcal{T}_w$ with full support can be constructed by setting $\mu_w:=\lambda^{[0,1]}\circ \hat{w}^{-1}$, where $\lambda^{[0,1]}$ is the usual one-dimensional Lebesgue measure on $[0,1]$. Furthermore, given a pair $(w,u)$, where $w\in\mathcal{W}$ and $u=(u_i)_{i\geq 1}\in[0,1]^\mathbb{N}$, we define a sequence of vertices $(\sigma_i)_{i\geq 1}$ by setting $\sigma_i=\hat{w}(u_i)$ for each $i$. This allows us to construct a sequence of subtrees $\mathcal{T}_{w,u}(k)$ of $\mathcal{T}_w$ as at (\ref{Tkdef}). Note that we will usually suppress the dependence on $w$ and $u$ from the notation for all of these objects when it is clear which excursion and sequence is being considered.

\begin{defn}\label{gammadef} The set $\Gamma$ is the collection of pairs $(w,u)\in \mathcal{W}\times [0,1]^\mathbb{N}$ such that $\mu$ satisfies (\ref{mulower}) for some $\kappa>0$, the sequence $(u_i)_{i\geq 1}$ is dense in $[0,1]$, and the vertices $(\sigma_i)_{i\geq 1}$ are a dense collection of leaves of $\mathcal{T}$, distinct and not equal to $\rho$ for any $i$.
\end{defn}

\section{$\ell^1$ convergence}\label{ell1}

For $x=(x(1),x(2),\dots)\in\ell^1$, define $[[0,x]]_{sp}$ as in \cite{Aldous3} to be the union of line segments connecting $0$ to $(x(1),0,0,\dots)$ to $(x(1),x(2),0,0,\dots),\dots$. Fix $k\geq 1$ and suppose we are given distinct $x^{(1)},\dots,x^{(k)}\in\ell^1\backslash\{0\}$. Write $\mathbf{x}=(x^{(1)},\dots,x^{(k)})$ and set
\begin{equation}\label{Tx}
\mathcal{T}^{\mathbf{x}}:=\bigcup_{i=1}^k [[0,x^{(i)}]]_{sp},
\end{equation}
which is a compact real tree. We assume that every $x\in\{0,x^{(1)},\dots,x^{(k)}\}$ is a leaf of $\mathcal{T}^\mathbf{x}$. Define $\lambda^{\mathbf{x}}$ to be one-dimensional Hausdorff measure on $\mathcal{T}^{\mathbf{x}}$; in this section it will be convenient not to normalise this to be a probability measure. Denote the law of the Brownian motion $X^{\mathbf{x}}$ on $(\mathcal{T}^{\mathbf{x}},\lambda^{\mathbf{x}})$ started from 0 by $\mathbf{P}_0^{\mathbf{x}}$. Given a Borel probability measure $\nu$ on $\mathcal{T}^\mathbf{x}$, let $A^{\mathbf{x},\nu}$ be the additive functional defined by
\[A^{\mathbf{x},\nu}_t:=\int_{\mathcal{T}^\mathbf{x}}L^\mathbf{x}_t(x)\nu(dx),\]
where $(L^\mathbf{x}_t(x))_{x\in\mathcal{T}^{\mathbf{x}},t\geq 0}$ are the local times of $X^\mathbf{x}$, which exist and are jointly continuous, $\mathbf{P}_0^\mathbf{x}$-a.s., because $\lambda^{\mathbf{x}}$ satisfies (\ref{mulower}) for $\kappa=1$ and therefore Lemma \ref{ltcont} applies.

The aim of this section is to show that if we have a sequence $(\mathbf{x}_n)_{n\geq 1}$, where $\mathbf{x}_n=(x_n^{(1)},\dots,x_n^{(k)})\in(\ell^1)^k$, that satisfies $\mathbf{x}_n\rightarrow\mathbf{x}$, then the Brownian motions $X^{\mathbf{x}_n}$ on $(\mathcal{T}^{\mathbf{x}_n},\lambda^{\mathbf{x}_n})$ started from 0 converge in distribution to the Brownian motion $X^{\mathbf{x}}$ started from 0. Note that we construct $\mathcal{T}^{\mathbf{x}_n}$ from $\mathbf{x}_n$ similarly to the definition of $\mathcal{T}^{\mathbf{x}}$ at (\ref{Tx}) and $\lambda^{\mathbf{x}_n}$ is the one-dimensional Hausdorff measure on $\mathcal{T}^{\mathbf{x}_n}$; the assumption that $\mathbf{x}_n\rightarrow \mathbf{x}$ means that we can define the law $\mathbf{P}_0^{\mathbf{x}_n}$ of the Brownian motion on $(\mathcal{T}^{\mathbf{x}_n},\lambda^{\mathbf{x}_n})$ started from 0 as in the previous section, at least for large $n$. Moreover, simultaneously with this convergence, we will prove that if $(\nu_n)_{n\geq 1}$ is a sequence of Borel probability measures, where $\nu_n$ is supported on $\mathcal{T}^{\mathbf{x}_n}$, such that $\nu_n\rightarrow\nu$ weakly as probability measures on $\ell^1$, then the related additive functionals $A^{\mathbf{x}_n,\nu_n}$, defined by
\[A^{\mathbf{x}_n,\nu_n}_t:=\int_{\mathcal{T}^{\mathbf{x}_n}}L^{\mathbf{x}_n}_t(x)\nu_n(dx),\]
converge in distribution to $A^{\mathbf{x},\nu}$ in $C(\mathbb{R}_+,\mathbb{R}_+)$. To deduce the existence and continuity of the local times $(L^{\mathbf{x}_n}_t(x))_{x\in\mathcal{T}^{\mathbf{x}_n},t\geq 0}$ of $X^{\mathbf{x}_n}$, we again apply Lemma \ref{ltcont}.

On line segments, the processes $X^{\mathbf{x}_n}$ and $X^\mathbf{x}$ look like standard one-dimensional Brownian motion. However, the structure of $\mathcal{T}^{\mathbf{x}_n}$ and $\mathcal{T}^{\mathbf{x}}$ can vary at branch-points, and so we must be careful about analysing the processes close to these. In our arguments, we will approximate $X^{\mathbf{x}_n}$ and $X^\mathbf{x}$ by processes that avoid the branch-points of the trees. Define the finite set of `vertices' of $\mathcal{T}^{\mathbf{x}}$ by
\[\mathcal{B}^{\mathbf{x}}:=\{b^{\mathbf{x}}(x,x',x''):x,x',x''\in\{0,x^{(1)},\dots,x^{(k)}\}\},\]
where $b^\mathbf{x}(x,x',x'')$ is the branch-point of $x$, $x'$ and $x''$ in $\mathcal{T}^{\mathbf{x}}$. Note that $\mathcal{B}^\mathbf{x}$ contains the set $\{0,x^{(1)},\dots,x^{(k)}\}$. An $\varepsilon$-neighbourhood of $\mathcal{B}^{\mathbf{x}}$ in $\ell^1$ is given by
\[\mathcal{B}^{\mathbf{x}}_\varepsilon:=\bigcup_{x\in\mathcal{B}^{\mathbf{x}}}B_{\ell^1}(x,\varepsilon),\]
where $B_{\ell^1}(x,\varepsilon)$ is the open ball in $(\ell^1,d_{\ell^1})$ of radius $\varepsilon$ centred at $x$. Now, fix a strictly positive constant $\varepsilon_1<\varepsilon_0$, where
\begin{equation}\label{eps0}
\varepsilon_0:=\tfrac{1}{2}\inf_{\substack{x,x'\in\mathcal{B}^{\mathbf{x}}:\\x\neq x'}}d_{\ell^1}(x,x').
\end{equation}
Set $\upsilon_0=0$ and, for $i\geq 0$, let
\begin{eqnarray*}
\varsigma_i&:=&\inf\left\{t\geq \upsilon_{i}:X^\mathbf{x}_t\not\in\mathcal{B}^{\mathbf{x}}_{\varepsilon_1}\right\},\\
\upsilon_{i+1}&:=&\inf\left\{t\geq \varsigma_{i}:X^\mathbf{x}_t\in\overline{\mathcal{B}^{\mathbf{x}}_{\varepsilon_1/2}}\right\}.
\end{eqnarray*}
Define
\[A^{\mathbf{x},\varepsilon_1}_t:=\int_{0}^t\mathbf{1}_{\{s\in I\}}ds,\]
where
\begin{equation}\label{idef}
I:=\mathbb{R}\backslash \cup_{i=0}^\infty (\upsilon_i,\varsigma_i),
\end{equation}
and its inverse $\tau^{\mathbf{x},\varepsilon_1}(t):=\inf\{s: A^{\mathbf{x},\varepsilon_1}_s>t\}$.
Finally, let $X^{\mathbf{x},\varepsilon_1}$ be the process defined by
\[X^{\mathbf{x},\varepsilon_1}_t:=X^{\mathbf{x}}_{\tau^{\mathbf{x},\varepsilon_1}(t)},\]
which takes values in the space $D(\mathbb{R}_+,\ell^1)$ of cadlag paths in $\ell^1$. To prove that $X^{\mathbf{x},\varepsilon_1}$ approximates $X^{\mathbf{x}}$ well for small $\varepsilon_1$, we will apply the following (rather crude) bound on the expectations of the local times $(L^{\mathbf{x}}_t(x))_{x\in\mathcal{T}^\mathbf{x},t\geq0}$.

{\lem \label{local}The local times of $X^{\mathbf{x}}$ satisfy
\[\sup_{x\in\mathcal{T}^\mathbf{x}}\mathbf{E}_0^\mathbf{x}L^\mathbf{x}_t(x)\leq \frac{t+4\lambda^\mathbf{x}(\mathcal{T}^\mathbf{x})e^t}{\lambda^\mathbf{x}(\mathcal{T}^\mathbf{x})}.\]}
\begin{proof} By \cite{BlumGet}, V.3.28, we have $\mathbf{P}_0^\mathbf{x}\left(|L^\mathbf{x}_t(x)-L^\mathbf{x}_t(x')|>2\delta\right)\leq 2e^te^{-\delta}$, for every $x,x'\in\mathcal{T}^\mathbf{x}$. Integrating this inequality implies that $\mathbf{E}_0^\mathbf{x}\left(|L^\mathbf{x}_t(x)-L^\mathbf{x}_t(x')|\right)\leq 4e^t$. Thus
\begin{eqnarray*}
\mathbf{E}_0^\mathbf{x}\left(L^\mathbf{x}_t(x)\right)\lambda^\mathbf{x}(\mathcal{T}^\mathbf{x})&\leq& \mathbf{E}_0^\mathbf{x}\left(\int_{\mathcal{T}^\mathbf{x}} L^\mathbf{x}_t(x')\lambda^\mathbf{x}(dx')\right)\\
&&+
\int_{\mathcal{T}^\mathbf{x}} \mathbf{E}_0^\mathbf{x}\left(|L^\mathbf{x}_t(x)-L^\mathbf{x}_t(x')|\right)\lambda^\mathbf{x}(dx')\\
&\leq&t+4e^t \lambda^\mathbf{x}(\mathcal{T}^\mathbf{x}).
\end{eqnarray*}
\end{proof}

{\lem \label{xeps1} For $t_0,\varepsilon>0$,
\[\lim_{\varepsilon_1\rightarrow0}\mathbf{P}_0^\mathbf{x}\left(\sup_{t\in[0,t_0]}d_{\ell^1}(X^\mathbf{x}_t,X^{\mathbf{x},\varepsilon_1}_t)>\varepsilon\right)=0.\]}
\begin{proof} By construction $X^{\mathbf{x},\varepsilon_1}_t=X^{\mathbf{x}}_{\tau^{\mathbf{x},\varepsilon_1}(t)}$. Hence, applying the continuity of $X^\mathbf{x}$ and the definition of $\tau^{\mathbf{x},\varepsilon_1}$ as the inverse of $A^{\mathbf{x},\varepsilon_1}$, the result will follow if we can show that
\begin{equation}\label{atight1}
\lim_{\varepsilon_1\rightarrow0}\mathbf{P}_0^\mathbf{x}\left(\sup_{t\in[0,t_0+1]}\left|t-A^{\mathbf{x},\varepsilon_1}_t\right|>\varepsilon\right)=0.
\end{equation}
To prove this first note that if $X^{\mathbf{x}}_0=0$, then
\[\sup_{t\in[0,t_0+1]}\left|t-A^{\mathbf{x},\varepsilon_1}_t\right|\leq \int_{0}^{t_0+1}\mathbf{1}_{\{X_s\in\mathcal{B}^\mathbf{x}_{\varepsilon_1}\}}ds=\int_{\mathcal{B}^\mathbf{x}_{\varepsilon_1}}L^\mathbf{x}_{t_0+1}(x)\lambda^{\mathbf{x}}(dx).\]
Thus, applying Markov's inequality, Fubini's theorem and Lemma \ref{local}, it follows that
\begin{equation}\label{abound}
\mathbf{P}_0^\mathbf{x}\left(\sup_{t\in[0,t_0+1]}\left|t-A^{\mathbf{x},\varepsilon_1}_t\right|>\varepsilon\right)\leq \varepsilon^{-1}\int_{\mathcal{B}^\mathbf{x}_{\varepsilon_1}}\mathbf{E}_0^\mathbf{x}\left(L^\mathbf{x}_{t_0+1}(x)\right)\lambda^{\mathbf{x}}(dx)\leq c\lambda^{\mathbf{x}}\left(\mathcal{B}^\mathbf{x}_{\varepsilon_1}\right),
\end{equation}
where $c$ is a constant that does not depend on $\varepsilon_1$. Since $\lambda^\mathbf{x}$ is non-atomic and $\mathcal{B}^\mathbf{x}$ is a finite set, the result follows.
\end{proof}

We now prove a similar result for $X^{\mathbf{x}_n}$ that is uniform in $n$. Set
$\upsilon^n_0=0$ and, for $i\geq 0$, let
\begin{eqnarray}
\varsigma^n_i&:=&\inf\left\{t\geq \upsilon^n_{i}:X^{\mathbf{x}_n}_t\not\in\mathcal{B}^{\mathbf{x}}_{\varepsilon_1}\label{taun}\right\},\\
\upsilon^n_{i+1}&:=&\inf\left\{t\geq \varsigma^n_{i}:X^{\mathbf{x}_n}_t\in\overline{\mathcal{B}^{\mathbf{x}}_{\varepsilon_1/2}}\right\}.\label{upsn}
\end{eqnarray}
Note that although $\mathcal{B}^{\mathbf{x}}_{\varepsilon_1/2}$ and $\mathcal{B}^{\mathbf{x}}_{\varepsilon_1}$ do not depend on $n$, the above quantities will be well-defined and finite for large $n$; when they are not, simply set $\varsigma^n_0=0$ and $\varsigma^n_i=\upsilon_i^n=\infty$ for $i\geq 1$. Define $A^{\mathbf{x}_n,\varepsilon_1}$, $\tau^{\mathbf{x}_n,\varepsilon_1}$ and $X^{\mathbf{x}_n,\varepsilon_1}$ from these stopping times analogously to the definitions of $A^{\mathbf{x},\varepsilon_1}$, $\tau^{\mathbf{x},\varepsilon_1}$ and $X^{\mathbf{x},\varepsilon_1}$ respectively.

{\lem \label{xneps1}If $\mathbf{x}_n\rightarrow\mathbf{x}$, then, for $t_0,\varepsilon>0$,
\[\lim_{\varepsilon_1\rightarrow0}\limsup_{n\rightarrow\infty}\mathbf{P}_0^{\mathbf{x}_n}\left(\sup_{t\in[0,t_0]}d_{\ell^1}(X^{\mathbf{x}_n}_t,X^{\mathbf{x}_n,\varepsilon_1}_t)>\varepsilon\right)=0.\]}
\begin{proof} Similarly to the bound at (\ref{abound}), it is possible to deduce that
\begin{equation}\label{aaaaa}
\limsup_{n\rightarrow\infty}\mathbf{P}_0^{\mathbf{x}_n}\left(\sup_{t\in[0,t_0+1]}\left|t-A^{\mathbf{x}_n,\varepsilon_1}_t\right|>\varepsilon\right)\leq \limsup_{n\rightarrow\infty}c\lambda^{\mathbf{x}_n}\left(\mathcal{B}^\mathbf{x}_{\varepsilon_1}\right)= c\lambda^{\mathbf{x}}\left(\mathcal{B}^\mathbf{x}_{\varepsilon_1}\right),
\end{equation}
where $c$ is a constant that does not depend on $n$ or $\varepsilon_1$. The result will follow from this if we can show that the sequence $(\mathbf{P}_0^{\mathbf{x}_n})_{n\geq 1}$ is tight. To prove this, first observe that it is elementary to check that there exist constants $c_1,c_2\in(0,\infty)$ such that
\[c_1r\leq \liminf_{n\rightarrow\infty}\inf_{x\in\mathcal{T}^{\mathbf{x}_n}}\lambda^{\mathbf{x}_n}\left(B_{\ell^1}(x,r)\right)\leq \limsup_{n\rightarrow\infty}\sup_{x\in\mathcal{T}^{\mathbf{x}_n}}\lambda^{\mathbf{x}_n}\left(B_{\ell^1}(x,r)\right)\leq c_2r\]
for every $r\in(0,1]$. By applying the argument of \cite{Kumagai}, Lemma 4.2, this implies that
\[\limsup_{n\rightarrow\infty}\sup_{x\in\mathcal{T}^{\mathbf{x}_n}}\mathbf{P}_x^{\mathbf{x}_n}\left(\inf\{s:d_{\ell^1}\left(x,X^{\mathbf{x}_n}_s\right)>r\}<t\right)\leq c_3e^{-\frac{c_4r^2}{t}}\]
for every $r\in(0,1]$, $t\in(0,t_1]$, for some constants $c_3,c_4,t_1\in(0,\infty)$. Consequently
\[\lim_{t\rightarrow0}\limsup_{n\rightarrow\infty}t^{-1}\sup_{x\in\mathcal{T}^{\mathbf{x}_n}}\mathbf{P}_x^{\mathbf{x}_n}\left(\inf\{s:d_{\ell^1}\left(x,X^{\mathbf{x}_n}_s\right)>r\}<t\right)=0\]
for any $r>0$, which implies the tightness of $(\mathbf{P}_0^{\mathbf{x}_n})_{n\geq 1}$ as desired (cf. the corollary to Theorem 7.4 of \cite{Bill2}).
\end{proof}

We now construct an approximation for the additive functional $A^{\mathbf{x},\nu}$. First, formulate a local time $L^{\mathbf{x},\varepsilon_1}$ for the process $X^{\mathbf{x},\varepsilon_1}$ by setting
\begin{equation}\label{ltap}
L^{\mathbf{x},\varepsilon_1}_t(x):=\int_{0}^{\tau^{\mathbf{x},\varepsilon_1}(t)}\mathbf{1}_{\{s\in I\}}d_sL^{\mathbf{x}}_s(x)
\end{equation}
for $x\in\mathcal{T}^\mathbf{x}$ and $t\geq 0$, where $I$ is defined as at (\ref{idef}), then let $A^{\mathbf{x},\nu,\varepsilon_1}$ be defined by
\begin{equation}\label{axne}
A^{\mathbf{x},\nu,\varepsilon_1}_t:=\int_{\mathcal{T}_{\varepsilon_1}^\mathbf{x}}L^{\mathbf{x},\varepsilon_1}_t(x)\nu(dx)+\sum_{x\in\mathcal{B}^\mathbf{x}}\nu(B_{\ell_1}(x,\varepsilon_1))\sup_{y\in\partial B_{\ell^1}(x,\varepsilon_1)\cap\mathcal{T}^\mathbf{x}}L^{\mathbf{x},\varepsilon_1}_t(y),
\end{equation}
where $\mathcal{T}^\mathbf{x}_{\varepsilon_1}:=\mathcal{T}^{\mathbf{x}}\backslash \mathcal{B}^\mathbf{x}_{\varepsilon_1}$. That $A^{\mathbf{x},\nu,\varepsilon_1}$ is uniformly close to $A^{\mathbf{x},\nu}$ is confirmed by the following lemma.

{\lem \label{adfcl} For $t_0,\varepsilon>0$,
\[\lim_{\varepsilon_1\rightarrow0}\mathbf{P}_0^\mathbf{x}\left(\sup_{t\in[0,t_0]}|A_t^{\mathbf{x},\nu}-A_t^{\mathbf{x},\nu,\varepsilon_1}|>\varepsilon\right)=0.\]}
\begin{proof} Since $\{s:X_s\in\mathcal{T}^\mathbf{x}_{\varepsilon_1}\}\subseteq I$, the definition of $L^{\mathbf{x},\varepsilon_1}$ at (\ref{ltap}) implies that $L^{\mathbf{x},\varepsilon_1}_t(x)=L^\mathbf{x}_{\tau^{\mathbf{x},\varepsilon_1}(t)}(x)$ for every $x\in\mathcal{T}_{\varepsilon_1}^\mathbf{x}$ and $t\geq 0$. Consequently, one can check that
\[\sup_{t\in[0,t_0]}|A_t^{\mathbf{x},\nu}-A_t^{\mathbf{x},\nu,\varepsilon_1}|\leq\sup_{t\in[0,t_0]}\sup_{\substack {x,y\in\mathcal{T}^\mathbf{x}:\\d_{\ell^1}(x,y)\leq 2 \varepsilon_1}}|L^{\mathbf{x}}_t(x)-L^\mathbf{x}_{\tau^{\mathbf{x},\varepsilon_1}(t)}(y)|,\]
where we apply the fact that $\nu$ is a probability measure. The $\mathbf{P}_0^\mathbf{x}$-a.s. joint continuity of the local times $L^\mathbf{x}$ allows us to deduce the result from (\ref{atight1}).
\end{proof}

For large $n$, we obtain objects analogous to $L^{\mathbf{x},\varepsilon_1}$ and $A^{\mathbf{x},\nu,\varepsilon_1}$ by setting
\[L^{\mathbf{x}_n,\varepsilon_1}_t(x):=\int_{0}^{\tau^{\mathbf{x}_n,\varepsilon_1}(t)}\mathbf{1}_{\{s\in I_n\}}d_sL^{\mathbf{x}_n}_s(x)\]
for $x\in\mathcal{T}^{\mathbf{x}_n}$ and $t\geq 0$, where $I_n=\mathbb{R}\backslash \cup_{i=0}^\infty (\upsilon^n_i,\varsigma^n_i)$ is defined from the stopping times introduced at (\ref{taun}) and (\ref{upsn}), then letting $A^{\mathbf{x}_n,\nu_n,\varepsilon_1}$ be defined by
\begin{equation}\label{axnen}
A^{\mathbf{x}_n,\nu_n,\varepsilon_1}_t:=\int_{\mathcal{T}_{\varepsilon_1}^{\mathbf{x}_n}}L^{\mathbf{x}_n,\varepsilon_1}_t(x)\nu_n(dx)+\sum_{x\in\mathcal{B}^\mathbf{x}}\nu_n(B_{\ell_1}(x,\varepsilon_1))\sup_{y\in\partial B_{\ell^1}(x,\varepsilon_1)\cap\mathcal{T}^{\mathbf{x}_n}}L^{\mathbf{x}_n,\varepsilon_1}_t(y),
\end{equation}
where $\mathcal{T}^{\mathbf{x}_n}_{\varepsilon_1}:=\mathcal{T}^{\mathbf{x}_n}\backslash \mathcal{B}^\mathbf{x}_{\varepsilon_1}$ and the summands where $B_{\ell^1}(x,\varepsilon_1)\cap\mathcal{T}^{\mathbf{x}_n}=\emptyset$ are assumed to be equal to zero. To check that $A^{\mathbf{x}_n,\nu_n,\varepsilon_1}$ is close to $A^{\mathbf{x}_n,\nu_n}$ uniformly in $n$, we will apply the following tightness result for the local times $(L^{\mathbf{x}_n})_{n\geq 1}$.

{\lem\label{lntight} If $\mathbf{x}_n\rightarrow\mathbf{x}$, then, for $t_0,\varepsilon>0$,
\begin{equation}\label{prob}
\lim_{\delta\rightarrow0}\limsup_{n\rightarrow\infty}\mathbf{P}_0^{\mathbf{x}_n}\left( \sup_{\substack{x,y\in\mathcal{T}^{\mathbf{x}_n}:\\d_{\ell^1}(x,y)\leq \delta}}
\sup_{\substack{s,t\in[0,t_0]:\\|s-t|\leq \delta}} |L^{\mathbf{x}_n}_s(x)-L^{\mathbf{x}_n}_t(y)|>\varepsilon\right)=0.
\end{equation}}
\begin{proof} First note that any two vertices $x,y\in\mathcal{T}^{\mathbf{x}_n}$ are contained in a set of the form $[[0,x_n^{(i)}]]^{\mathbf{x}_n}$ or $[[x_n^{(i)},x_n^{(j)}]]^{\mathbf{x}_n}$ for some $i,j\in\{1,\dots,k\}$, $i\neq j$, where for $x',y'\in\mathcal{T}^{\mathbf{x}_n}$, we write $[[x',y']]^{\mathbf{x}_n}$ to represent the path from $x'$ to $y'$ in $\mathcal{T}^{\mathbf{x}_n}$. Secondly, writing $b^{\mathbf{x}_n}(x,x',x'')$ to represent the branch-point of $x$, $x'$ and $x''$ in $\mathcal{T}^{\mathbf{x}_n}$, applying the assumption $\mathbf{x}_n\rightarrow \mathbf{x}$ and equation (5) of \cite{Aldous3}, it is possible to check that there exists a finite integer $n_0$ such that
\[\eta(n):=\sup_{i,j=1,\dots,k} d_{\ell^1}\left(b^{\mathbf{x}_n}(0,x^{(i)}_n,x^{(j)}_n),b^\mathbf{x}(0,x^{(i)},x^{(j)})\right)< \varepsilon_0/2\]
for $n\geq n_0$, where $\varepsilon_0$ is defined at (\ref{eps0}). This implies that $d_{\ell^1}(0,x_n^{(i)})> \varepsilon_0$ and also $d_{\ell^1}(b^{\mathbf{x}_n}(0,x_n^{(i)},x_n^{(j)}),x_n^{(i)})> \varepsilon_0$ for every $i\neq j$, whenever $n\geq n_0$. For $n\geq n_0$, it is an elementary exercise to deduce from these two facts the existence of a collection of paths $([[a_i,b_i]]^{\mathbf{x}_n})_{i\in\mathcal{I}_n}$, where $d_{\ell^1}(a_i,b_i)=\varepsilon_0$ and $d_{\ell^1}(b^{\mathbf{x}_n}(0,a_i,b_i),a_i)\in\{0,\varepsilon_0/2,\varepsilon_0\}$, such that: if $x,y\in\mathcal{T}^{\mathbf{x}_n}$ and $d_{\ell^1}(x,y)<\varepsilon_0/2$, then $x$ and $y$ are both contained in a single set of the form $[[a_i,b_i]]^{\mathbf{x}_n}$ for some $i\in\mathcal{I}_n$. Moreover, the collection $([[a_i,b_i]]^{\mathbf{x}_n})_{i\in\mathcal{I}_n}$ can be chosen in such a way that $\#\mathcal{I}_n$ is uniformly bounded in $n$.

Suppose now that $n\geq n_0$ is large enough so that a cover of the form described above exists and $\delta<\varepsilon_0/2$. For $i\in\mathcal{I}_n$, define
\[A^{i}_t:=\int_{[[a_i,b_i]]^{\mathbf{x}_n}}L^{\mathbf{x}_n}_t(x)\lambda^{\mathbf{x}_n}(dx),\]
and $\tau^i(t):=\inf\{s:A^i_s>t\}$. Similarly to Lemma \ref{timechange}, the law of the process $X^i$, where $X^i_t:=X^{\mathbf{x}_n}_{\tau^i(t)}$, under $\mathbf{P}_0^{\mathbf{x}_n}$ is equal to the law of the Brownian motion on the measured compact real tree $([[a_i,b_i]]^{\mathbf{x}_n},\lambda^{\mathbf{x}_n}([[a_i,b_i]]^{\mathbf{x}_n}\cap\cdot))$ started from $b^{\mathbf{x}_n}(0,a_i,b_i)$. Moreover, the local times of $X^i$ are given by $L^i_t(x)=L^{\mathbf{x}_n}_{\tau^i(t)}(x)$ for $x\in[[a_i,b_i]]^{\mathbf{x}_n}$ and $t\geq 0$ (cf. \cite{Croydoncbp}, Lemma 3.4). Now, by our choice of the intervals $([[a_i,b_i]]^{\mathbf{x}_n})_{i\in\mathcal{I}_n}$, we have that
\begin{eqnarray*}
\lefteqn{\sup_{\substack{x,y\in\mathcal{T}^{\mathbf{x}_n}:\\d_{\ell^1}(x,y)\leq \delta}}
\sup_{\substack{s,t\in[0,t_0]:\\|s-t|\leq \delta}} |L^{\mathbf{x}_n}_s(x)-L^{\mathbf{x}_n}_t(y)|}\\
&\leq&\sum_{i\in\mathcal{I}_n}\sup_{\substack{x,y\in [[a_i,b_i]]^{\mathbf{x}_n}:\\d_{\ell^1}(x,y)\leq \delta}}
\sup_{\substack{s,t\in[0,t_0]\cap\overline{\tau^i(\mathbb{R}_+)}:\\|s-t|\leq \delta}} |L^{\mathbf{x}_n}_s(x)-L^{\mathbf{x}_n}_t(y)|,
\end{eqnarray*}
where the condition $s,t\in \overline{\tau^i(\mathbb{R}_+)}$ is justified by the observation that, for $x\in[[a_i,b_i]]^{\mathbf{x}_n}$, the local time $L^{\mathbf{x}_n}_t(x)$ only increases on the set $\{t\geq 0:X_t=x\}\subseteq \overline{\tau^i(\mathbb{R}_+)}$. Furthermore, simple continuity arguments allow us to replace $\overline{\tau^i(\mathbb{R}_+)}$ by ${\tau^i(\mathbb{R}_+)}$. Since $|A^i_s-A^i_t|\leq |s-t|$, we have that $|\tau^i(s)-\tau^i(t)|\geq |s-t|$, and consequently we obtain from this bound that
\[\sup_{\substack{x,y\in\mathcal{T}^{\mathbf{x}_n}:\\d_{\ell^1}(x,y)\leq \delta}}
\sup_{\substack{s,t\in[0,t_0]:\\|s-t|\leq \delta}} |L^{\mathbf{x}_n}_s(x)-L^{\mathbf{x}_n}_t(y)|\leq
\sum_{i\in\mathcal{I}_n}\sup_{\substack{x,y\in [[a_i,b_i]]^{\mathbf{x}_n}:\\d_{\ell^1}(x,y)\leq \delta}}
\sup_{\substack{s,t\in[0,t_0]:\\|s-t|\leq \delta}} |L^{i}_s(x)-L^{i}_t(y)|,\]
where we apply the representation of $(L^i_t(x))_{x\in[[a_i,b_i]]^{\mathbf{x}_n},t\geq0}$ noted above. All the measured real trees $\{([[a_i,b_i]]^{\mathbf{x}_n},\lambda^{\mathbf{x}_n}([[a_i,b_i]]^{\mathbf{x}_n}\cap\cdot))\}_{i\in\mathcal{I}_n}$ are equivalent and, by construction, the Brownian motion $X^{\mathbf{x}_n}$ first hits $[[a_i,b_i]]^{\mathbf{x}_n}$ at $b^{\mathbf{x}_n}(0,a_i,b_i)$; hence it follows that the probability in the left-hand side of (\ref{prob}) is bounded above by
\[\#\mathcal{I}_n \sup_{x\in\{0,\varepsilon_0/2\}}\mathbf{P}_x^{[0,\varepsilon_0]}\left(\sup_{\substack{y,z\in[0,\varepsilon_0]:\\|y-z|\leq \delta}}
\sup_{\substack{s,t\in[0,t_0]:\\|s-t|\leq \delta}} |L^{[0,\varepsilon_0]}_s(y)-L^{[0,\varepsilon_0]}_t(z)|>\varepsilon(\#\mathcal{I}_n)^{-1}\right),\]
where $\mathbf{P}_x^{[0,\varepsilon_0]}$ is the law of the Brownian motion on the real interval $[0,\varepsilon_0]$ equipped with the one-dimensional Hausdorff measure and $(L^{[0,\varepsilon_0]}_t(x))_{x\in[0,\varepsilon_0],t\geq 0}$ are the local times of this process. Recalling that $\#\mathcal{I}_n $ is uniformly bounded in $n$, Lemma \ref{ltcont} allows it to be deduced that the above expression decays to zero as first $n\rightarrow\infty$ and then $\delta\rightarrow 0$.
\end{proof}

By following a similar proof to that of Lemma \ref{adfcl}, combining the previous lemma with (\ref{aaaaa}) allows it to be to deduced that $A^{\mathbf{x}_n,\nu_n,\varepsilon_1}$ does indeed approximate $A^{\mathbf{x}_n,\nu_n}$ as desired.

{\lem \label{adfcln} If $\mathbf{x}_n\rightarrow\mathbf{x}$, then, for $t_0,\varepsilon>0$,
\[\lim_{\varepsilon_1\rightarrow0}\limsup_{n\rightarrow \infty}\mathbf{P}_0^{\mathbf{x}_n}\left(\sup_{t\in[0,t_0]}|A_t^{\mathbf{x}_n,\nu_n}-A_t^{\mathbf{x}_n,\nu_n,\varepsilon_1}|>\varepsilon\right)=0.\]}
\bigskip

\begin{figure}[t]
\centering
\scalebox{0.4}{\includegraphics{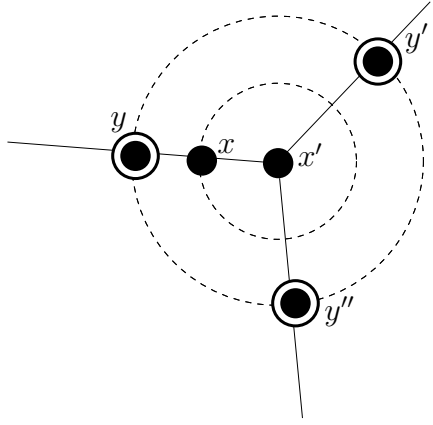}}
\put(-80,100){$x$}
\put(-50,95){$x'$}
\put(-120,110){$y$}
\put(-10,140){$y'$}
\put(-40,37){$y''$}
\caption{Example of $x\in \partial\mathcal{T}^\mathbf{x}_{\varepsilon_1/2}$, $x'\in\mathcal{B}^\mathbf{x}$, $N^\mathbf{x}_{\varepsilon_1}(x)=\{y,y',y''\}$.}
\label{bp}
\end{figure}

Our next step is to demonstrate that $X^{\mathbf{x}_n,\varepsilon_1}$ converges to $X^{\mathbf{x},\varepsilon_1}$ and $A^{\mathbf{x}_n,\nu_n,\varepsilon_1}$ is close to $A^{\mathbf{x},\nu,\varepsilon_1}$, which we do by applying a particular sample path construction of the processes from their excursions between jump-times. Let us start by introducing some further notation. Set $\mathcal{T}^\mathbf{x}_{\varepsilon_1/2}:=\mathcal{T}^{\mathbf{x}}\backslash \mathcal{B}^\mathbf{x}_{\varepsilon_1/2}$, which consists of a finite number of connected components, each a closed line segment with end-points in the finite set $\partial\mathcal{T}^\mathbf{x}_{\varepsilon_1/2}$. For $x\in\partial\mathcal{T}^\mathbf{x}_{\varepsilon_1/2}$, write
\[N^\mathbf{x}_{\varepsilon_1}(x):=\{y\in\partial\mathcal{T}_{\varepsilon_1}^\mathbf{x}:[[x,y]]^\mathbf{x}\cap\partial\mathcal{T}_{\varepsilon_1}^\mathbf{x}=\{y\}\}\]
to represent the collection of `nearest neighbours' of $x$ in $\partial\mathcal{T}_{\varepsilon_1}^\mathbf{x}$, where $[[x,y]]^\mathbf{x}$ is the path from $x$ to $y$ in $\mathcal{T}^\mathbf{x}$. For a pictorial representation of the definitions at a typical branch-point, see Figure \ref{bp}. For $x\in\mathcal{T}^\mathbf{x}_{\varepsilon_1/2}$, we will denote by $C^\mathbf{x}_{\varepsilon_1}(x)$ the connected component of $\mathcal{T}_{\varepsilon_1/2}^\mathbf{x}$ containing $x$.

Let $(\zeta_i)_{i\geq 0}$ be the jump-times of $X^{\mathbf{x},\varepsilon_1}$, by convention we set $\zeta_0=0$. Applying the definition of Brownian motion on a dendrite, conditional on $X_{\zeta_i}^{\mathbf{x},\varepsilon_1}$, the path segment $(X^{\mathbf{x},\varepsilon_1}_{\zeta_i+t})_{t\in[0,\zeta_{i+1}-\zeta_i)}$ is distributed precisely as a standard one-dimensional Brownian motion on the line segment $C_{\varepsilon_1}^\mathbf{x}(X^{\mathbf{x},\varepsilon_1}_{\zeta_i})$, started from $X^{\mathbf{x},\varepsilon_1}_{\zeta_i}$ and run until it hits $\partial\mathcal{T}^\mathbf{x}_{\varepsilon_1/2}$. The local times $(L^{\mathbf{x},\varepsilon_1}_{\zeta_i+t}(x))_{x\in C_{\varepsilon_1}^\mathbf{x}(X^{\mathbf{x},\varepsilon_1}_{\zeta_i}),t\in[0,\zeta_{i+1}-\zeta_i)}$ are distributed exactly as the local times of the same one-dimensional Brownian motion; outside $C_{\varepsilon_1}^\mathbf{x}(X^{\mathbf{x},\varepsilon_1}_{\zeta_i})$, the local times $L^{\mathbf{x},\varepsilon_1}$ do not increase on the time interval $[\zeta_i,\zeta_{i+1})$. Moreover, the strong Markov property implies that, conditional on $X_{\zeta_i}^{\mathbf{x},\varepsilon_1}$, $(X^{\mathbf{x},\varepsilon_1}_{\zeta_i+t},(L^{\mathbf{x},\varepsilon_1}_{\zeta^i+t}(x))_{x\in\mathcal{T}^{\mathbf{x}}})_{t\in[0,\zeta_{i+1}-\zeta_i)}$ is independent of the $\sigma$-algebra generated by $(X^{\mathbf{x},\varepsilon_1}_t)_{t\leq\zeta_i}$. At discontinuities, the process satisfies the following transition law:
\[p^{\mathbf{x},\varepsilon_1}(x,y):=\mathbf{P}^\mathbf{x}_0\left(X^{\mathbf{x},\varepsilon_1}_{\zeta_i}=y|X_{\zeta_i^{-}}^{\mathbf{x},\varepsilon_1}=x\right)=\left\{\begin{array}{ll}
             \frac{1+\#N^\mathbf{x}_{\varepsilon_1}(x)}{2\#N^\mathbf{x}_{\varepsilon_1}(x)}, & \mbox{if $y\in C_{\varepsilon_1}^{\mathbf{x}}(x)$,}\\
             \frac{1}{2\#N^\mathbf{x}_{\varepsilon_1}(x)}, & \mbox{otherwise,}
           \end{array}\right.\]
for $x\in \partial\mathcal{T}^\mathbf{x}_{\varepsilon_1/2}$, $y\in N^\mathbf{x}_{\varepsilon_1}(x)$. As a consequence of this description, we can construct $(X^{\mathbf{x},\varepsilon_1},L^{\mathbf{x},\varepsilon_1})$ from a countable collection
\begin{equation}\label{collection}
\left\{\left(\alpha^{x,i}\right)_{x\in\partial\mathcal{T}^\mathbf{x}_{\varepsilon_1/2}},\left(\beta^{y,i},\gamma^{y,i},\xi^{y,i}\right)_{y\in \partial\mathcal{T}^\mathbf{x}_{\varepsilon_1}}\right\}_{i\geq 0}
\end{equation}
of random variables built on an underlying probability space with probability measure $\mathbf{P}$ that satisfy the following properties:
\begin{itemize}
  \item The random variables $\left(\alpha^{x,i}\right)_{x\in\partial\mathcal{T}^\mathbf{x}_{\varepsilon_1/2},i\geq 0}$ are independent. The random variable $\alpha^{x,i}$ is $N^\mathbf{x}_{\varepsilon_1}(x)$-valued and distributed according to the law determined by $p^{\mathbf{x},\varepsilon_1}(x,\cdot)$.
  \item The triples $\left\{\left(\beta^{y,i},\gamma^{y,i},\xi^{y,i}\right)\right\}_{y\in \partial\mathcal{T}^\mathbf{x}_{\varepsilon_1},i\geq 0}$ are independent of each other and of the collection $\left(\alpha^{x,i}\right)_{x\in\partial\mathcal{T}^\mathbf{x}_{\varepsilon_1/2},i\geq 0}$. The process $\beta^{y,i}=(\beta^{y,i}_t)_{t\geq 0}$ is a Brownian motion on the line segment $C^\mathbf{x}_{\varepsilon_1}(y)$ (equipped with the appropriate restriction of $\lambda^{\mathbf{x}}$) started from $y$; $\gamma^{y,i}:=\inf\{s:\beta^{y,i}_s\in \partial\mathcal{T}^\mathbf{x}_{\varepsilon_1/2}\}$; and $\xi^{y,i}=(\xi^{y,i}_t(x))_{x\in C^\mathbf{x}_{\varepsilon_1}(y),t\geq 0}$ are the jointly continuous local times of $\beta^{y,i}$.
\end{itemize}
We now define a cadlag process $\tilde{X}^{\mathbf{x},\varepsilon_1}$, its times of discontinuities $(\tilde{\zeta}_i)_{i\geq1}$ and its local times $\tilde{L}^{\mathbf{x},\varepsilon_1}$ from these random variables. First, set $\tilde{X}^{\mathbf{x},\varepsilon_1}_0$ to be the unique point in $\mathcal{T}^\mathbf{x}$ satisfying $d_{\ell^1}(0,\tilde{X}^{\mathbf{x},\varepsilon_1}_0)=\varepsilon_1$ (such a point is well-defined, because we are assuming that 0 is a leaf of $\mathcal{T}^\mathbf{x}$ and $\varepsilon_1$ is less than $\varepsilon_0$, as defined at (\ref{eps0})), $\tilde{\zeta}_0=0$ and $\tilde{L}^{\mathbf{x},\varepsilon_1}_0(x)=0$ for every $x\in\mathcal{T}^\mathbf{x}$. Given $(\tilde{X}^{\mathbf{x},\varepsilon_1}_t,(\tilde{L}^{\mathbf{x},\varepsilon_1}_t(x))_{x\in\mathcal{T}^\mathbf{x}})_{t\in[0,\tilde{\zeta}_i]}$, define
\[\tilde{\zeta}_{i+1}:=\tilde{\zeta}_i+\gamma^{\tilde{X}^{\mathbf{x},\varepsilon_1}_{\tilde{\zeta}_i},i},\]
\[\left(\tilde{X}^{\mathbf{x},\varepsilon_1}_{t}\right)_{t\in(\tilde{\zeta_i},\tilde{\zeta}_{i+1})}:=
\left(\beta^{\tilde{X}^{\mathbf{x},\varepsilon_1}_{\tilde{\zeta}_i},i}_{t-\tilde{\zeta}_i}\right)_{t\in(\tilde{\zeta_i},\tilde{\zeta}_{i+1})},\hspace{20pt}\tilde{X}^{\mathbf{x},\varepsilon_1}_{\tilde{\zeta}_{i+1}}:=\alpha^{\tilde{X}^{\mathbf{x},\varepsilon_1}_{\tilde{\zeta}_{i+1}^-},i},\]
and also
\[\left(\tilde{L}^{\mathbf{x},\varepsilon_1}_t(x)\right)_{t\in(\tilde{\zeta_i},\tilde{\zeta}_{i+1}]}:=\left\{\begin{array}{ll}
                                                                                                   \left(\tilde{L}^{\mathbf{x},\varepsilon_1}_{\tilde{\zeta}_i}(x)+\xi^{\tilde{X}^{\mathbf{x},\varepsilon_1}_{\tilde{\zeta}_i},i}_{t-\tilde{\zeta}_i}(x)\right)_{t \in(\tilde{\zeta_i},\tilde{\zeta}_{i+1}]},&\mbox{if $x\in C_{\varepsilon_1}^\mathbf{x}(\tilde{X}^{\mathbf{x},\varepsilon_1}_{\tilde{\zeta}_i})$}, \\
                                                                                                   \left(\tilde{L}^{\mathbf{x},\varepsilon_1}_{\tilde{\zeta}_i}(x)\right)_{t\in(\tilde{\zeta_i},\tilde{\zeta}_{i+1}]},&\mbox{otherwise.}
                                                                                                \end{array}\right.\]
Applying the above description of $(X^{\mathbf{x},\varepsilon_1},L^{\mathbf{x},\varepsilon_1})$, it is straightforward to check that the law of $(\tilde{X}^{\mathbf{x},\varepsilon_1},\tilde{L}^{\mathbf{x},\varepsilon_1})$ under $\mathbf{P}$ is identical to the law of $(X^{\mathbf{x},\varepsilon_1},L^{\mathbf{x},\varepsilon_1})$ under $\mathbf{P}_0^\mathbf{x}$.

\begin{figure}[t]
\centering
\scalebox{0.8}{\includegraphics{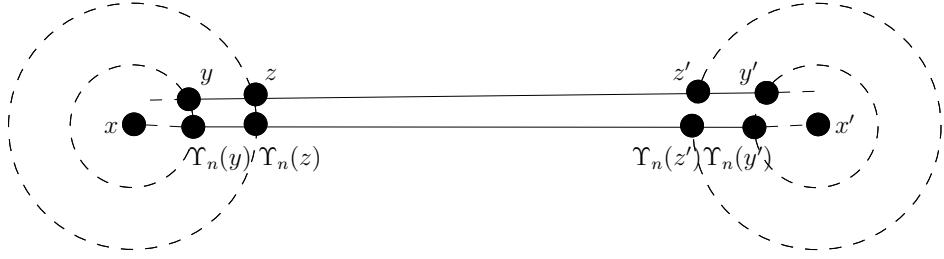}
\put(-50,55){$x'$}
\put(-95,80){$y'$}
\put(-112,40){$\Upsilon_n(y')$}
\put(-126,80){$z'$}
\put(-145,40){$\Upsilon_n(z')$}
\put(-392,55){$x$}
\put(-347,80){$y$}
\put(-353,40){$\Upsilon_n(y)$}
\put(-317,80){$z$}
\put(-320,40){$\Upsilon_n(z)$}
}
\caption{Example of homeomorphism between line-segments.}
\label{edge}
\end{figure}

We continue by describing the structure of $\mathcal{T}^{\mathbf{x}_n}_{\varepsilon_1/2}$. Define $\eta(n)$ as in the proof of Lemma \ref{lntight}. The assumption $\mathbf{x}_n\rightarrow\mathbf{x}$ implies that there exists a finite integer $n_0$ such that $\eta(n)<\varepsilon_1/2$ for $n\geq n_0$, which implies in turn that $\mathcal{T}^{\mathbf{x}_n}_{\varepsilon_1/2}$ is homeomorphic to $\mathcal{T}^\mathbf{x}_{\varepsilon_1/2}$ for $n\geq n_0$. In particular, if we suppose that $x$ and $x'$ are neighbours in $\mathcal{B}^{\mathbf{x}}$, by which we mean that $x,x'\in\mathcal{B}^\mathbf{x}$, $x\neq x'$ and $[[x,x']]^\mathbf{x}\cap\mathcal{B}^\mathbf{x}=\{x,x'\}$, then there exists a unique connected component of $\mathcal{T}^{\mathbf{x}_n}_{\varepsilon_1/2}$ which is a closed line segment with end-points in $\partial B_{\ell^1}(x,\varepsilon_1/2)$ and
$\partial B_{\ell^1}(x',\varepsilon_1/2)$; moreover, every connected component of $\mathcal{T}^{\mathbf{x}_n}_{\varepsilon_1/2}$ can be represented in this way. We define a homeomorphism from $\mathcal{T}^{\mathbf{x}_n}_{\varepsilon_1/2}$ to $\mathcal{T}^{\mathbf{x}}_{\varepsilon_1/2}$ that: maps the end-point of such a line segment, $y$ say, in $\partial B_{\ell^1}(x,\varepsilon_1/2)$ to the unique point $\Upsilon_n(y)$ in $[[x,x']]^\mathbf{x}$ that satisfies $d_{\ell^1}(x,\Upsilon_n(y))=\varepsilon_1/2$; maps the point of the line segment at a distance $\varepsilon_1/2$ from $y$, $z$ say, to the unique point $\Upsilon_n(z)$ in $[[x,x']]^\mathbf{x}$ that satisfies $d_{\ell^1}(x,\Upsilon_n(z))=\varepsilon_1$; and $\Upsilon_n$ is extended by linear interpolation on the line-segments between points for which we have not already defined it. Figure \ref{edge} depicts a typical configuration on a line-segment.

We now assume that $\eta(n)<\varepsilon_1/2$. The transition law at the discontinuities $(\zeta^n_i)_{i\geq 1}$ of $X^{\mathbf{x}_n,\varepsilon_1}$ is given by
\[p^{\mathbf{x}_n,\varepsilon_1}(x,y):=\mathbf{P}^{\mathbf{x}_n}_0\left(X^{\mathbf{x}_n,\varepsilon_1}_{\zeta^n_i}=y|X^{\mathbf{x}_n,\varepsilon_1}_{{\zeta^n_i}^{-}}=x\right),\]
where $x\in \partial \mathcal{T}^{\mathbf{x}_n}_{\varepsilon_1/2}$ and $y\in N_{\varepsilon_1}^{\mathbf{x}_n}(x):= \Upsilon_n^{-1}(N^\mathbf{x}_{\varepsilon_1}(\Upsilon_n(x)))$. We use this to define a cadlag process $\tilde{X}^{\mathbf{x}_n,\varepsilon_1}$, for each $n$, from a countable collection
\[\left\{\left(\alpha^{n,x,i}\right)_{x\in\partial\mathcal{T}^{\mathbf{x}_n}_{\varepsilon_1/2}},\left(\beta^{n,y,i},\gamma^{n,y,i},\xi^{n,y,i}\right)_{y\in \partial\mathcal{T}^{\mathbf{x}_n}_{\varepsilon_1}}\right\}_{i\geq 0}\]
of random variables built on our underlying probability space with probability measure $\mathbf{P}$ that satisfy the following properties, where we apply the notation $C^{\mathbf{x}_n}_{\varepsilon_1}(x)$ to represent the connected component of $\mathcal{T}^{\mathbf{x}_n}_{\varepsilon_1/2}$ containing $x\in\mathcal{T}^{\mathbf{x}_n}_{\varepsilon_1/2}$:
\begin{itemize}
  \item The random variables $\left(\alpha^{n,x,i}\right)_{x\in\partial\mathcal{T}^{\mathbf{x}_n}_{\varepsilon_1/2},i\geq 0}$ are independent. The random variable $\alpha^{n,x,i}$ is $N^{\mathbf{x}_n}_{\varepsilon_1}(x)$-valued and distributed according to $p^{\mathbf{x}_n,\varepsilon_1}(x,\cdot)$.
  \item The triples $\left\{\left(\beta^{n,y,i},\gamma^{n,y,i},\xi^{n,y,i}\right)\right\}_{y\in \partial\mathcal{T}^{\mathbf{x}_n}_{\varepsilon_1},i\geq 0}$ are independent of each other and of the collection $\left(\alpha^{n,x,i}\right)_{x\in\partial\mathcal{T}^{\mathbf{x}_n}_{\varepsilon_1/2},i\geq 0}$. The process $\beta^{n,y,i}=(\beta^{n,y,i}_t)_{t\geq 0}$ is a Brownian motion on the line segment $C_{\varepsilon_1}^{\mathbf{x}_n}(y)$ (equipped with the appropriate restriction of $\lambda^{\mathbf{x}_n}$) started from $y$; $\gamma^{n,y,i}:=\inf\{s:\beta^{n,y,i}_s\in \partial\mathcal{T}^{\mathbf{x}_n}_{\varepsilon_1/2}\}$; and $\xi^{n,y,i}=(\xi^{n,y,i}_t(x))_{x\in C^{\mathbf{x}_n}_{\varepsilon_1}(y),t\geq 0}$ are the jointly continuous local times of $\beta^{n,y,i}$.
\end{itemize}
Constructing $(\tilde{X}^{\mathbf{x}_n,\varepsilon_1},\tilde{L}^{\mathbf{x}_n,\varepsilon_1})$ from these random variables, similarly to the definition of $(\tilde{X}^{\mathbf{x},\varepsilon_1},\tilde{L}^{\mathbf{x},\varepsilon_1})$, results in a process whose law under $\mathbf{P}$ is equal to that of $({X}^{\mathbf{x}_n,\varepsilon_1},{L}^{\mathbf{x}_n,\varepsilon_1})$ under $\mathbf{P}_0^{\mathbf{x}_n}$. The following simple lemma is crucial in proving that ${X}^{\mathbf{x}_n,\varepsilon_1}$ converges to ${X}^{\mathbf{x},\varepsilon_1}$ and $A^{\mathbf{x}_n,\nu_n,\varepsilon_1}$ is close to $A^{\mathbf{x},\nu,\varepsilon_1}$ as we do in the subsequent result. Note that, the topology on the countable collections with the same index sets we consider is that generated by uniform convergence of finite sub-collections (with respect to the appropriate product topology).

{\lem \label{family}If $\mathbf{x}_n\rightarrow\mathbf{x}$, the countable collection of random variables consisting of
\[\left(\Upsilon_n(\alpha^{n,\Upsilon_n^{-1}(x),i})\right)_{x\in\partial\mathcal{T}^\mathbf{x}_{\varepsilon_1/2},i\geq0}\]
and
\[\left(\beta^{n,\Upsilon_n^{-1}(y),i},\gamma^{n,\Upsilon_n^{-1}(y),i},\left(\xi^{n,\Upsilon_n^{-1}(y),i}(\Upsilon_n^{-1}(x))\right)_{x\in C_{\varepsilon_1}^\mathbf{x}(y)}\right)_{y\in \partial\mathcal{T}^\mathbf{x}_{\varepsilon_1},i\geq 0}\]
converges in distribution to the collection at (\ref{collection}).}
\begin{proof} To deduce that $\Upsilon_n(\alpha^{n,\Upsilon_n^{-1}(x),i})$ converges in distribution to $\alpha^{x,i}$, it suffices to check that
\[p^{\mathbf{x}_n,\varepsilon_1}(\Upsilon_n^{-1}(x),\Upsilon_n^{-1}(y))\rightarrow p^{\mathbf{x},\varepsilon_1}(x,y)\]
for each $x\in \partial\mathcal{T}^\mathbf{x}_{\varepsilon_1/2}$, $y\in N_{\varepsilon_1}^{\mathbf{x}}(x)$. The proof that this is true involves an elementary electrical network (or harmonic analysis) calculation and is omitted (that the Brownian motion on a real tree can be constructed in terms of electrical resistance networks is guaranteed by \cite{Croydoncbp}, Proposition 2.2, and a detailed study of harmonic analysis on such spaces appears in \cite{Kigamidendrite}).

The convergence of the triple
\[\left(\beta^{n,\Upsilon_n^{-1}(y),i},\xi^{n,\Upsilon_n^{-1}(y),i},(\gamma^{n,\Upsilon_n^{-1}(y),i}(\Upsilon_n^{-1}(x)))_{x\in C_{\varepsilon_1}^\mathbf{x}(y)}\right)\]
to $(\beta^{y,i},\gamma^{y,i},\xi^{y,i})$ is a simple consequence of the fact that
\[\left((B_{t})_{t\geq 0}, h(-a_n,b_n),(L_{t}(x))_{x\in\mathbb{R},t\geq 0}\right)\rightarrow\left((B_{t})_{t\geq 0},h(-a,b),(L_{t}(x))_{x\in\mathbb{R},t\geq 0}\right)\]
in distribution  whenever $a_n\rightarrow a$ and $b_n\rightarrow b$, for some $a,b>0$, where: $B$ is a standard Brownian motion in $\mathbb{R}$ started from 0; $h(-a,b)$ is the hitting time of $\{-a,b\}$ by $B$; and $L$ are the jointly continuous local times of $B$. Note that to map this result into our setting, we apply the fact that $\sup_{x\in\mathcal{T}^{\mathbf{x}_n}}d_{\ell^1}(\Upsilon_n(x),x)\rightarrow 0$ as $n\rightarrow\infty$.
\end{proof}

{\lem \label{approxi} Suppose $\mathbf{x}_n\rightarrow\mathbf{x}$. If $F$ is a continuous bounded function on $D(\mathbb{R}_+,\ell^1)$, then
\begin{equation}\label{first}
\lim_{n\rightarrow\infty}\mathbf{E}_0^{\mathbf{x}_n}\left(F(X^{\mathbf{x}_n,\varepsilon_1})\right)=\mathbf{E}_0^{\mathbf{x}}\left(F(X^{\mathbf{x},\varepsilon_1})\right).
\end{equation}
Moreover, if $\nu_n\rightarrow \nu$ weakly as probability measures on $\ell^1$ and $F$ is a continuous bounded function on $D(\mathbb{R}_+,\ell^1)\times C(\mathbb{R}_+,\mathbb{R}_+)$, then
\begin{equation}\label{second}
\lim_{\varepsilon_1\rightarrow 0}\limsup_{n\rightarrow\infty}\left|\mathbf{E}_0^{\mathbf{x}_n}\left(F(X^{\mathbf{x}_n,\varepsilon_1},A^{\mathbf{x}_n,\nu_n,\varepsilon_1})\right)-\mathbf{E}_0^{\mathbf{x}}\left(F(X^{\mathbf{x},\varepsilon_1},A^{\mathbf{x},\nu,\varepsilon_1})\right)\right|=0.
\end{equation}}
\begin{proof} By separability, it is possible to assume that all the relevant random variables are constructed in such a way that the convergence of Lemma \ref{family} holds $\mathbf{P}$-a.s. It will follow easily from this that $\tilde{X}^{\mathbf{x}_n,\varepsilon_1}\rightarrow\tilde{X}^{\mathbf{x},\varepsilon_1}$ in $D(\mathbb{R}_+,\ell^1)$, $\mathbf{P}$-a.s., if we can show that the number of discontinuities that $\tilde{X}^{\mathbf{x},\varepsilon_1}$ admits in any finite time interval is finite. By continuity, there exists an $\varepsilon>0$ such that $\mathbf{P}(\inf_{y\in\partial\mathcal{T}^{\mathbf{x}}_{\varepsilon_1}}\gamma^{y,i}>\varepsilon)>0$. Thus the strong law of large numbers implies that
\begin{equation}\label{findisc}
\lim_{i\rightarrow\infty}\tilde{\zeta}_i\geq\sum_{i=0}^{\infty}\inf_{y\in\partial\mathcal{T}^{\mathbf{x}}_{\varepsilon_1}}\gamma^{y,i}=\infty,
\end{equation}
$\mathbf{P}$-a.s., and hence there can indeed only be a finite number of discontinuities of $\tilde{X}^{\mathbf{x},\varepsilon_1}$ in any finite time interval. Due to the equivalence of the laws of $\tilde{X}^{\mathbf{x}_n,\varepsilon_1},\tilde{X}^{\mathbf{x},\varepsilon_1}$ and ${X}^{\mathbf{x}_n,\varepsilon_1},{X}^{\mathbf{x},\varepsilon_1}$ under the appropriate probability measures, this yields the convergence result at (\ref{first}).

Let us now suppose we have a realisation of random variables such that the convergence of Lemma \ref{family} occurs and (\ref{findisc}) holds. Under these conditions, it is possible to check that, for any $t_0>0$,
\begin{equation}\label{epsdec}
\varepsilon(n):=\sup_{t\in[0,t_0]}\sup_{x\in\mathcal{T}^{\mathbf{x}}_{\varepsilon_1}}\left|\tilde{L}^{\mathbf{x}_n,\varepsilon_1}_t(\Upsilon_n^{-1}(x))-\tilde{L}^{\mathbf{x},\varepsilon_1}_t(x)\right|\rightarrow0,
\end{equation}
as $n\rightarrow \infty$. Moreover, we can define a function $(\bar{L}^{\mathbf{x},\varepsilon_1}_t(x))_{x\in\mathcal{T}^\mathbf{x}\cup\mathcal{B}^\mathbf{x}_{\varepsilon_1},t\geq 0}$ that is jointly continuous and agrees with $\tilde{L}^{\mathbf{x},\varepsilon_1}_t(x)$ for $x\in\mathcal{T}^{\mathbf{x}}_{\varepsilon_1}$. If $x\in B_{\ell^1}(y,\varepsilon_1)$ for some $y\in\mathcal{B}^\mathbf{x}$, then we set
\[\bar{L}_t^{\mathbf{x},\varepsilon_1}(x):=\sum_{z\in\partial B_{\ell^1}(y,\varepsilon_1)\cap\mathcal{T}^\mathbf{x}}d_{\ell^1}(x,z)^{-1}\tilde{L}^{\mathbf{x},\varepsilon_1}_t(z)/\sum_{z\in\partial B_{\ell^1}(y,\varepsilon_1)\cap\mathcal{T}^\mathbf{x}}d_{\ell^1}(x,z)^{-1},\]
so that $\bar{L}_t^{\mathbf{x},\varepsilon_1}(x)$ is a weighted average of the points in $\mathcal{T}^{\mathbf{x}}$ on the boundary of $B_{\ell^1}(y,\varepsilon_1)$. Furthermore, for large $n$, let $\phi_n:\mathcal{T}^{\mathbf{x}_n}\cup\mathcal{B}^\mathbf{x}_{\varepsilon_1}\rightarrow \mathcal{T}^{\mathbf{x}}\cup\mathcal{B}^\mathbf{x}_{\varepsilon_1}$ be defined by setting $\phi_n(x)=\Upsilon_n(x)$ on $\mathcal{T}^{\mathbf{x}_n}_{\varepsilon_1}$ and $\phi_n(x)=x$ otherwise, and note that $\nu_n\circ\phi_n^{-1}\rightarrow\nu$ weakly as probability measures on $\mathcal{T}^{\mathbf{x}}\cup\mathcal{B}^\mathbf{x}_{\varepsilon_1}$. If we construct $\tilde{A}^{\mathbf{x}_n,\nu_n,\varepsilon_1}$ and $\tilde{A}^{\mathbf{x},\nu,\varepsilon_1}$ from $\tilde{L}^{\mathbf{x}_n,\varepsilon_1}$ and $\tilde{L}^{\mathbf{x},\varepsilon_1}$ analogously to the definitions of ${A}^{\mathbf{x}_n,\nu_n,\varepsilon_1}$ and ${A}^{\mathbf{x},\nu,\varepsilon_1}$ at (\ref{axne}) and (\ref{axnen}) respectively, then we can apply the above notation to deduce that, for any $0\leq t_1\leq \dots\leq t_m\leq t_0$,
\begin{eqnarray*}
\lefteqn{\sup_{t\in \{t_1,\dots,t_m\}}\left|\tilde{A}_t^{\mathbf{x}_n,\nu_n,\varepsilon_1}-\tilde{A}_t^{\mathbf{x},\nu,\varepsilon_1}\right|}\\
&\leq& \varepsilon(n)+2\sup_{t\in[0,t_0]}\sup_{\substack{x,y\in\mathcal{T}^\mathbf{x}\cup\mathcal{B}^\mathbf{x}_{\varepsilon_1}:\\d_{\ell^1}(x,y)\leq 2\varepsilon_1}}\left|\bar{L}_t^{\mathbf{x},\varepsilon_1}(x)-\bar{L}_t^{\mathbf{x},\varepsilon_1}(y)\right|\\
&& +\sup_{t\in \{t_1,\dots,t_m\}}\left|\int_{\mathcal{T}^\mathbf{x}\cup\mathcal{B}_{\varepsilon_1}^\mathbf{x}}\bar{L}_t^{\mathbf{x},\varepsilon_1}(x)\nu_n\circ\phi_n^{-1}(dx)-\int_{\mathcal{T}^\mathbf{x}\cup\mathcal{B}_{\varepsilon_1}^\mathbf{x}}\bar{L}_t^{\mathbf{x},\varepsilon_1}(x)\nu(dx)\right|
\end{eqnarray*}
By (\ref{epsdec}), the first term in the upper bound decays to zero as $n\rightarrow\infty$. The third term also converges to zero, since $\nu_n\circ\phi_n^{-1}\rightarrow\nu$. It follows from this and the definition of $\bar{L}^{\mathbf{x},\varepsilon_1}$ that, for any $\varepsilon>0$,
\begin{eqnarray*}
\lefteqn{\mathbf{P}\left(\limsup_{n\rightarrow\infty}\sup_{t\in \{t_1,\dots,t_m\}}\left|\tilde{A}_t^{\mathbf{x}_n,\nu_n,\varepsilon_1}-\tilde{A}_t^{\mathbf{x},\nu,\varepsilon_1}\right|>\varepsilon\right)}\\
&\leq & \mathbf{P}_0^\mathbf{x}\left(2\sup_{t\in[0,t_0]}\sup_{\substack{x,y\in\mathcal{T}^\mathbf{x}_{\varepsilon_1}:\\d_{\ell^1}(x,y)\leq 4\varepsilon_1}}\left|{L}_t^{\mathbf{x},\varepsilon_1}(x)-{L}_t^{\mathbf{x},\varepsilon_1}(y)\right|>\varepsilon\right)\\
&=&\mathbf{P}_0^\mathbf{x}\left(2\sup_{t\in[0,t_0]}\sup_{\substack{x,y\in\mathcal{T}^\mathbf{x}_{\varepsilon_1}:\\d_{\ell^1}(x,y)\leq 4\varepsilon_1}}\left|{L}_{\tau^{\mathbf{x},\varepsilon_1}(t)}^{\mathbf{x}}(x)-{L}_{\tau^{\mathbf{x},\varepsilon_1}(t)}^{\mathbf{x}}(y)\right|>\varepsilon\right),
\end{eqnarray*}
where the equality follows from the observation made in the proof of Lemma \ref{adfcl} that $L^{\mathbf{x},\varepsilon_1}_t(x)=L^\mathbf{x}_{\tau^{\mathbf{x},\varepsilon_1}(t)}(x)$ for $x\in\mathcal{T}^\mathbf{x}_{\varepsilon_1}$. By applying the $\mathbf{P}_0^{\mathbf{x}}$-a.s. continuity of the local times $L^{\mathbf{x}}$ and (\ref{atight1}), this implies that
\[\lim_{\varepsilon_1\rightarrow0}\mathbf{P}\left(\limsup_{n\rightarrow\infty}\sup_{t\in \{t_1,\dots,t_m\}}\left|\tilde{A}_t^{\mathbf{x}_n,\nu_n,\varepsilon_1}-\tilde{A}_t^{\mathbf{x},\nu,\varepsilon_1}\right|>\varepsilon\right)=0.\]
In conjunction with the result of the opening paragraph of the proof and Lemma \ref{adfcln}, this will imply the convergence at (\ref{second}) once the tightness of $(A^{\mathbf{x}_n,\nu_n})_{n\geq 1}$ in $C(\mathbb{R}_+,\mathbb{R}_+)$ has been checked. However, this is an easy corollary of Lemma \ref{lntight}, and so the proof is complete.
\end{proof}

The main result of this section, which we can now state precisely, is an immediate consequence of Lemmas \ref{xeps1}, \ref{xneps1}, \ref{adfcl}, \ref{adfcln} and \ref{approxi}.

{\prop\label{main3} Suppose $\mathbf{x}_n\rightarrow\mathbf{x}$ and $\nu_n\rightarrow \nu$ weakly as probability measures on $\ell^1$. If $F$ is a continuous bounded function on $C(\mathbb{R}_+,\ell^1)\times C(\mathbb{R}_+,\mathbb{R}_+)$, then
\[\lim_{n\rightarrow\infty}\mathbf{E}_0^{\mathbf{x}_n}\left(F(X^{\mathbf{x}_n},A^{\mathbf{x}_n,\nu_n})\right)=\mathbf{E}_0^{\mathbf{x}}\left(F(X^{\mathbf{x}},A^{\mathbf{x},\nu})\right).\]}

\section{Simple random walk convergence}\label{discconv}

The results of the two previous sections allow us to prove our main scaling limit theorems by relatively straightforward adaptations of the proofs in \cite{Croydoncbp}. Since most of the objects we study here were also considered in \cite{Croydoncbp}, we will be brief in introducing them, and refer the reader to \cite{Croydoncbp} for further details.

Let $(T_n)_{n\geq 1}$ be a collection of ordered graph trees such that $T_n$ has $n$ vertices for each $n$. The root of $T_n$ will be denoted $\rho=\rho(T_n)$. Define the depth-first search $\hat{w}_n:\{0,\dots,2n\}\rightarrow T_n$ as in \cite{Aldous3}, Section 2.6 ($\hat{w}_n$ is extended from the definition there by setting $\hat{w}_n(0)=\hat{w}_n(2n)=\rho$), and suppose the search-depth process $(w_n(t))_{t\in[0,1]}$ is the function satisfying
\[w_n(i/2n):=d_{{T}_n}(\rho,\hat{w}_n(i)),\hspace{20pt}0\leq i\leq 2n,\]
where $d_{{T}_n}$ is the graph distance on ${T}_n$, and which is linear between these values, so that $w_n$ takes values in $C([0,1],\mathbb{R}_+)$. The uniform probability measure on the vertices of $T_n$ will be written as $\mu_n$. For each $n$, if we construct a function $\gamma_n:[0,1]\rightarrow[0,1]$ by setting
\[\gamma_n(t):=\left\{\begin{array}{ll}
                        \lfloor 2nt \rfloor/2n, & \mbox{if }w_n(\lfloor 2nt \rfloor/2n)\geq w_n( \lceil 2nt \rceil/2n),\\
                        \lceil 2nt \rceil/2n, & \mbox{otherwise,}
                      \end{array}
\right.\]
then it is the case that
\begin{equation}\label{mundef}
\mu_n=\lambda^{[0,1]}\circ(2n\gamma_n)^{-1}\circ \hat{w}_n^{-1},
\end{equation}
where, as previously, $\lambda^{[0,1]}$ is the one-dimensional Lebesgue measure on $[0,1]$ (this result can be checked by arguing along the lines of \cite{Aldous3}, Lemma 12). It will useful for later to note that $\gamma_n$ satisfies $\sup_{t\in[0,1]}|\gamma_n(t)-t|\leq (2n)^{-1}$

For a sequence $(u^n_k)_{k\geq1}\in[0,1]^\mathbb{N}$, define $(\sigma^n_k)_{k\geq1}\in T_n^\mathbb{N}$ by $\sigma_k^n:=\hat{w}_n(2n\gamma_n(u_k^n))$. Set $T_n(k)$ to be the minimal graph tree spanning $\{\rho,\sigma_1^n,\dots,\sigma_k^n\}$. The measure projection of $\mu_n$ onto ${T}_n(k)$ is denoted
\[\mu_n^{(k)}:=\mu_n\circ\phi_{{T}_n,{T}_n(k)}^{-1},\]
where the projection operator $\phi_{{T}_n,{T}_n(k)}$ is defined on graph trees analogously to the projection operator for real trees (see Section \ref{bmrt}).

We will repeatedly apply the following assumption throughout the remainder of the section. The set $\Gamma$ was introduced in Definition \ref{gammadef}.

\begin{assu}\label{a1} For each $n$, the sequence $(u^n_k)_{k\geq 1}$ is dense in $[0,1]$. Furthermore, there exists a divergent sequence $(\alpha_n)_{n\geq1}$ such that $\alpha_n=o(n)$ and
\[\left(\alpha_n^{-1}w_n,u^n\right)\rightarrow (w,u),\]
in $C([0,1],\mathbb{R}_+)\times [0,1]^\mathbb{N}$, for some $(w, u)\in \Gamma$.
\end{assu}

As in Section \ref{bmrt} for real trees, we use the sequential construction of \cite{Aldous3}, Section 2.2, to isometrically embed the vertices of ${T}_n$ into $\ell^1$ from the vertex sequence $(\sigma^n_k)_{k\geq 1}$. Observe that, under Assumption \ref{a1}, because we are assuming $(u^n_k)_{k\geq 1}$ to be dense in $[0,1]$, the sequence $(\sigma^n_k)_{k\geq 1}$ will contain all the vertices of ${T}_n$, and so this procedure does result in an isometric embedding for $T_n$. We shall denote by $\psi_n$ the unique distance-preserving map from the vertices of ${T}_n$ into the set $\{(x(1),x(2),\dots)\in\ell^1:x(i)\geq 0, i=1,2,\dots\}\subseteq\ell^1$ that satisfies $\psi_n(\rho)=0$ and
\begin{equation}\label{projn}
\pi_k(\psi_n(\sigma))=\psi_n(\phi_{{T}_n,{T}_n(k)}(\sigma))
\end{equation}
for every $\sigma\in T_n$ and $k\geq 1$, where $\pi_k$ is the projection operator defined below (\ref{proj}). We write $\tilde{{T}}_n, \tilde{\mu}_n,\dots$ to represent the $\ell^1$-embedded versions of objects.

Let us now introduce the discrete time processes that we will consider. Suppose that $X^n=(X^n_m)_{m\geq 0}$ is the discrete time simple random walk on $T_n$ started from $\rho$, and denote its law by $\mathbf{P}^{T_n}_\rho$. Set $X^{n,k}=\phi_{T_n,T_n(k)}(X^n)$, and let $J^{n,k}$ be the associated jump-chain, so that $J^{n,k}$ is the simple random walk on the vertices of $T_n(k)$ started from $\rho$. If $(A_m^{n,k})_{m\geq0}$ are the jump times of $X^{n,k}$, i.e. $A^{n,k}_0=0$ and, for $m\geq 1$,
\[A_m^{n,k}:=\min\left\{l\geq A_{m-1}^{n,k}:\:X^n_{l}\in{T}_n(k)\backslash\{X^{n}_{A_{m-1}^{n,k}}\}\right\},\]
and $(\tau^{n,k}(m))_{m\geq 0}$ is the discrete inverse of $A^{n,k}$, i.e. ${\tau}^{n,k}(m):=\max\{l:\:{A}^{n,k}_{l}\leq m\}$, then we can recover $X^{n,k}$ from $J^{n,k}$ through the identity
\begin{equation}\label{xnk}
X^{n,k}_m=J^{n,k}_{\tau^{n,k}(m)}.
\end{equation}
The local times of $J^{n,k}$ are determined by
\[L^{n,k}_m(\sigma):=\frac{2}{\mathrm{deg}_{n,k}(\sigma)}
\sum_{l=0}^{m}\mathbf{1}_\sigma(J_{l}^{n,k}),\]
for  $\sigma$ a vertex in ${T}_n(k)$, where $\mathrm{deg}_{n,k}(\sigma)$ is the usual graph degree of $\sigma$ in $T_n(k)$, and we use these to define an additive functional, $(\hat{A}^{n,k}_m)_{m\geq 0}$, by setting $\hat{A}^{n,k}_0=0$, and for $m\geq1$,
\[\hat{A}^{n,k}_m:=n\int_{{T}_n(k)}L^{n,k}_{m-1}(\sigma)\mu_n^{(k)}(d\sigma).\]
The discrete time inverse of $\hat{A}^{n,k}$ is given by $\hat{\tau}^{n,k}(m):=\max\{l:\hat{A}^{n,k}_l\leq m\}$, and we use this to define a time-changed version of $J^{n,k}$, denoted $(\hat{X}^{n,k}_m)_{m\geq 0}$, by setting
\begin{equation}\label{hatx}
\hat{X}^{n,k}_m:=J^{n,k}_{\hat{\tau}^{n,k}(m)}.
\end{equation}
As in \cite{Croydoncbp}, to show $\hat{X}^{n,k}$ and $X^{n,k}$ are close, we will prove a tightness result for $A^{n,k}$ and $\hat{A}^{n,k}$. We start by proving a continuity result for the local times $L^{n,k}$. Coinciding with the notation of \cite{Croydoncbp}, we define $\Lambda_n^{(k)}:=\alpha_n^{-1}\#E(T_n(k))$, where $E(T_n(k))$ is the edge set of $T_n(k)$.

{\lem\label{lnktight} Fix $k\in\mathbb{N}$ and $t_0>0$. Under Assumption \ref{a1},
\[\lim_{\delta\rightarrow0}\limsup_{n\rightarrow\infty}\mathbf{P}_\rho^{T_n}\left(\alpha_n^{-1} \sup_{\substack{\sigma,\sigma'\in T_n(k):\\d_{T_n}(\sigma,\sigma')\leq \delta\alpha_n}} \sup_{m\leq t_0\alpha_n^2\Lambda_n^{(k)}}\left|L^{n,k}_m(\sigma)-L^{n,k}_m(\sigma')\right|>\varepsilon\right)=0.\]}
\begin{proof} By assumption, $\alpha_n^{-1}d_{T_n}(\rho,\sigma^n_i)\rightarrow d_{\mathcal{T}}(\rho,\sigma_i)>0$ for each $i\geq 1$. Hence there exists a constant $L>0$ such that
\begin{equation}\label{Llower}
\inf_{i=1,\dots,k}\alpha_n^{-1}d_{T_n}(\rho,\sigma^n_i)\geq 2L
\end{equation}
for large $n$. For the remainder of the proof, we will suppose that $n$ is large enough so that this bound holds. If $\sigma\in T_n(k)$, then $\sigma$ is on the path between $\rho$ and $\sigma^n_i$ for some $i\in\{1,\dots,k\}$. Consequently, (\ref{Llower}) implies that there exists an injective path in $T_n(k)$ of length at least $\lfloor \alpha_n L\rfloor$ with endpoint $\sigma$. By considering the random walk observed on this path and the other vertices adjacent to $\sigma$, we can deduce that, for $t>0$,
\[\mathbf{P}_\rho^{T_n}\left(\alpha_n^{-1}L^{n,k}_{t_0\alpha_n^{2}\Lambda_n^{(k)}}(\sigma)\geq t\right)\leq \mathbf{P}_0^{\Gamma(\lfloor \alpha_n L\rfloor, \mathrm{deg}_{n,k}(\sigma)-1)}\left(\alpha_n^{-1}L^{\Gamma(\lfloor \alpha_n L\rfloor, \mathrm{deg}_{n,k}(\sigma)-1)}_{t_0\alpha_n^{2}\Lambda_n^{(k)}}\geq t\right),\]
where $\Gamma(i,D)$ is a graph tree consisting of a path of length $i$ emanating from a vertex $0$ along with $D$ other vertices each attached to $0$ by a single edge, and
$(L^{\Gamma(i,D)}_{m})_{m\geq 0}$ is the local time at $0$ of the simple random walk on $\Gamma(i,D)$. Using a strong Markov argument, it is possible to check that $L_m^{\Gamma(i,D)}$ is stochastically dominated by $\sum_{j=1}^{L_m^{\Gamma(i,0)}}\xi_j^D$, where $(\xi_j^D)_{j\geq 1}$ are independent geometric, parameter $D/(D+1)$, random variables, independent of the random walk on $\Gamma(i,0)$. Therefore
\begin{eqnarray}
{\mathbf{P}_\rho^{T_n}\left(\alpha_n^{-1}L^{n,k}_{t_0\alpha_n^{2}\Lambda_n^{(k)}}(\sigma)\geq t\right)}&\leq &\sup_{D\leq k}\mathbf{P}_0^{\Gamma(\lfloor \alpha_n L\rfloor, D)}\left(\alpha_n^{-1}L^{\Gamma(\lfloor \alpha_n L\rfloor, D)}_{t_0\alpha_n^{2}\Lambda_n^{(k)}}\geq t\right)\nonumber\\
&\leq &(t\alpha_n)^{-4} \sup_{D\leq k} \mathbf{E}_0^{\Gamma(\lfloor \alpha_n L\rfloor, 0)}\left(\sum_{i_1,i_2,i_3,i_4=1}^{L^{\Gamma(\lfloor \alpha_n L\rfloor, 0)}_{t_0\alpha_n^{2}\Lambda_n^{(k)}}}\xi_{i_1}^D \xi_{i_2}^D\xi_{i_3}^D\xi_{i_4}^D\right)\nonumber\\
&\leq& c_1 (t\alpha_n)^{-4} \mathbf{E}_0^{\Gamma(\lfloor \alpha_n L\rfloor, 0)}\left(\left(L^{\Gamma(\lfloor \alpha_n L\rfloor, 0)}_{t_0\alpha_n^{2}\Lambda_n^{(k)}}\right)^4\right)\nonumber\\
&\leq & c_2 t^{-4},\label{polybound}
\end{eqnarray}
where $c_1$ and $c_2$ are constants that do not depend on $\sigma$, $n$ or $t$. The final inequality here is an application of \cite{Croydoncbp}, Lemma B.2, and the easily checked fact that $\Lambda_n^{(k)}$ is uniformly bounded in $n$.

Given the above bound, it is possible to deduce that, for fixed $\varepsilon>0$,
\begin{equation}\label{locallocal}
\sup_{\substack{\sigma,\sigma'\in T_n(k):\\d_{T_n}(\sigma,\sigma')\leq \delta\alpha_n}} \mathbf{P}_\rho^{T_n}\left(\alpha_n^{-1} \sup_{m\leq t_0\alpha_n^2\Lambda_n^{(k)}}\left|L^{n,k}_m(\sigma)-L^{n,k}_m(\sigma')\right|>\varepsilon\right)\leq c_3\delta^2,
\end{equation}
for every $n\geq n_0$ and $\delta\in(0,1)$, where $n_0$ is a suitably large integer, by following the same argument as used to prove \cite{Croydoncbp}, Lemma 4.5 (which itself is an adaptation of a result appearing in \cite{Borodin}), inserting $\alpha_n$ in place of the scaling factor $n^{1/2}$ where appropriate. More specifically, fix $\sigma\neq\sigma'\in T_n(k)$ such that
$d_{T_n}(\sigma,\sigma')\leq \delta\alpha_n$ (note that in what follows we may assume that $\delta\alpha_n\geq 1$, else the left-hand side of (\ref{locallocal}) is clearly 0). Conditional on the event where the jump chain ${J}^{n,k}$ hits $\sigma$ before $\sigma'$ occurring, we have by a simple calculation
\begin{eqnarray}
\lefteqn{\sup_{m\leq t_0\alpha_n^2\Lambda_n^{(k)}} \left|{L}^{n,k}_{m}(\sigma)-{L}_m^{n,k}(\sigma')+2\sum_{i=1}^{\mathrm{deg}_{n,k}(\sigma){L}^{n,k}_{m}(\sigma)/2}
\eta_i\right|}\nonumber\\
&\leq&\sup
\left\{ 2N_i\mathrm{deg}_{n,k}(\sigma')^{-1} : i\leq \mathrm{deg}_{n,k}(\sigma){L}^{n,k}_{\lfloor t_0\alpha_n^2\Lambda_n^{(k)}\rfloor}(\sigma)/2\right\},\label{azbound}
\end{eqnarray}
where $N_i$ is the number of visits by ${J}^{n,k}$ to $\sigma'$ between the $i$th and $(i+1)$st visits to $\sigma$ and $\eta_i:=N_i\mathrm{deg}_{n,k}(\sigma')^{-1}-\mathrm{deg}_{n,k}(\sigma)^{-1}$ is a centred random variable (for the precise distribution of $N_i$ and estimates of the moments of $\eta_i$, see \cite{Croydoncbp}, Section B.2). Now, since $(\sum_{i=1}^m\eta_i)_{m\geq 1}$ is a martingale with respect to the filtration $(\mathcal{F}_m)_{m\geq 1}$, where $\mathcal{F}_m$ is
the $\sigma$-algebra generated by $J^{n,k}$ up to the $(m+1)$st hitting time of $\sigma$, and $L:=\mathrm{deg}_{n,k}(\sigma){L}^{n,k}_{\lfloor t_0\alpha_n^2\Lambda_n^{(k)}\rfloor}(\sigma)/2$ is a stopping time for this martingale, Doob's martingale norm inequality implies that
\begin{eqnarray*}
\mathbf{P}_\rho^{T_n}\left(2\alpha_n^{-1}\sup_{m\leq L}\left|\sum_{i=1}^{m}
\eta_i\right|>\varepsilon\right)&\leq&c_4\alpha_n^{-4}
\mathbf{E}_\rho^{T_n}\left(\left(\sum_{i=1}^{L}
\eta_i\right)^4\right).
\end{eqnarray*}
An upper bound for the right-hand side in terms of the moments of $L$ and $\eta_i$ can be computed by following the steps that lead to \cite{Borodin}, (1.29), which is an analogous bound for simple random walk on the line. In particular, one can check from (\ref{polybound}) and the estimates for the moments of $\eta_i$ of the form $|\mathbf{E}_\rho^{T_n}(\eta^p_i)|\leq c (\delta\alpha_n)^{p-1}$ appearing in the appendix of \cite{Croydoncbp} that
\begin{equation}\label{doobbound}
\mathbf{P}_\rho^{T_n}\left(2\alpha_n^{-1}\sup_{m\leq L}\left|\sum_{i=1}^{m}
\eta_i\right|>\varepsilon\right)\leq c_5\delta^2,
\end{equation}
uniformly in $n$, $\sigma$ and $\sigma'$. Moreover, it also holds that
\begin{eqnarray}
\mathbf{P}_\rho^{T_n}\left(\alpha_n^{-1}\sup_{i\leq L}2N_i\mathrm{deg}_{n,k}(\sigma')^{-1}>\varepsilon\right)&\leq&\mathbf{E}_\rho^{T_n}\left(L\right)\mathbf{P}_\rho^{T_n}\left(2\alpha_n^{-1}N_1\mathrm{deg}_{n,k}(\sigma')^{-1}
>\varepsilon\right)\nonumber\\
&\leq &c_6 \alpha_n^{-3}\left(1+\mathbf{E}_\rho^{T_n}\left(|\eta_1|^4\right)\right)\nonumber\\
&\leq & c_7 \delta^3,\label{otherbound}
\end{eqnarray}
uniformly in  $n$, $\sigma$ or $\sigma'$. Piecing together (\ref{azbound}), (\ref{doobbound}) and (\ref{otherbound}), we obtain the estimate (\ref{locallocal}) as desired.

Subsequently, by considering for each $\sigma\in T_n(k)$ the behaviour of the local times on the paths from $\sigma$ to the at most $k+1$ leaves of $T_n(k)$, we can apply a standard maximal inequality (for example, the extension of \cite{Bill2}, Theorem 10.3, suggested as \cite{Bill2}, Problem 10.1) to deduce from this result that
\[\sup_{\substack{\sigma\in T_n(k)}} \mathbf{P}_\rho^{T_n}\left(\alpha_n^{-1}\sup_{\substack{\sigma'\in T_n(k):\\d_{T_n}(\sigma,\sigma')\leq \delta\alpha_n}}  \sup_{m\leq t_0\alpha_n^2\Lambda_n^{(k)}}\left|L^{n,k}_m(\sigma)-L^{n,k}_m(\sigma')\right|>\varepsilon\right)\leq c_8\delta^2,\]
for every $n\geq n_0$ and $\delta\in(0,1)$ (cf. \cite{Croydoncbp}, equation (39)).

Finally, note that for each $n$ and $\delta$ we can choose a set $A_n(\delta)$ of at most $c_9\delta^{-1}$ vertices of $T_n(k)$ (where $c_9$ is independent of $n$ and $\delta$) such that the sets $\{\sigma':d_{T_n}(\sigma,\sigma')\leq \delta \alpha_n\}$, $\sigma\in A_n(\delta)$, cover $T_n(k)$. Therefore
\begin{eqnarray*}
\lefteqn{\mathbf{P}_\rho^{T_n}\left(\alpha_n^{-1} \sup_{\substack{\sigma,\sigma'\in T_n(k):\\d_{T_n}(\sigma,\sigma')\leq \delta\alpha_n}} \sup_{m\leq t_0\alpha_n^2\Lambda_n^{(k)}}\left|L^{n,k}_m(\sigma)-L^{n,k}_m(\sigma')\right|>\varepsilon\right)}\\
&\leq&\sum_{\sigma\in A_n(\delta)}2 \mathbf{P}_\rho^{T_n}\left(\alpha_n^{-1}\sup_{\substack{\sigma'\in T_n(k):\\d_{T_n}(\sigma,\sigma')\leq 2\delta\alpha_n}}  \sup_{m\leq t_0\alpha_n^2\Lambda_n^{(k)}}\left|L^{n,k}_m(\sigma)-L^{n,k}_m(\sigma')\right|>\varepsilon/2\right)\\
&\leq &c_{10} \delta,
\end{eqnarray*}
(cf. \cite{Croydoncbp}, Lemma 4.6), from which the lemma follows.
\end{proof}

We now show that the rescaled jump-chains $\tilde{J}^{n,k}:=\psi_n(J^{n,k})$ and additive functionals $\hat{A}^{n,k}$ converge. Henceforth, we extend $\tilde{J}^{n,k}$ and $\hat{A}^{n,k}$ to continuous time processes by linear interpolation in $\ell^1$ and $\mathbb{R}_+$ respectively. We recall the definition of $\mathcal{T}(k)$ from (\ref{Tkdef}) and also that the process $\tilde{X}^{\mathcal{T}(k),\lambda^{(k)}}$ is the $\ell^1$-embedding of ${X}^{\mathcal{T}(k),\lambda^{(k)}}$ (as introduced below Proposition \ref{prockconv}). The definition of $\hat{A}^{(k)}$ appears at (\ref{hatak}).

{\lem \label{jumpconv} Under Assumption \ref{a1}, the joint laws of the pairs
\[\left(\alpha_n^{-1}\tilde{J}^{n,k}_{t\alpha_n^2\Lambda_n^{(k)}}, (n\alpha_n)^{-1}\hat{A}^{n,k}_{t\alpha_n^{2}\Lambda_n^{(k)}}\right)_{t\geq 0}\]
under $\mathbf{P}_\rho^{T_n}$ converge to the joint law of
\[\left(\tilde{X}^{\mathcal{T}(k),\lambda^{(k)}}_t,\hat{A}^{(k)}_t\right)_{t\geq 0}\]
under $\mathbf{P}_\rho^{\mathcal{T}(k),\lambda^{(k)}}$ weakly as probability measures on the space $C(\mathbb{R}_+,\ell^1)\times C(\mathbb{R}_+,\mathbb{R}_+)$.}
\begin{proof} Fix $k\geq 1$. Consider the vectors $\mathbf{x}_n:=(\alpha_n^{-1}\psi_n(\sigma_1^n),\dots,\alpha_n^{-1}\psi_n(\sigma_k^n))$ and $\mathbf{x}:=(\psi(\sigma_1),\dots,\psi(\sigma_k))$. Under Assumption \ref{a1}, it is possible to check that $\mathbf{x}_n\rightarrow \mathbf{x}$. We can rewrite this as
\[\left(\alpha_n^{-1}\psi_n\left(\hat{w}_n\left(2n\gamma_n(u_i^n)\right)\right)\right)_{i\geq 1}\rightarrow \left(\psi\left(\hat{w}(u_i)\right)\right)_{i\geq 1}\]
as sequences in $\ell^1$. Since $(u_i)_{i\geq 1}$ is dense in $[0,1]$ and $\psi\circ \hat{w}$ is continuous, if we can show that $(\alpha_n^{-1}\psi_n\circ\hat{w}_n\circ(2n\gamma_n))_{n\geq 1}$ is tight, then we will obtain that
\begin{equation}\label{funcconv}
\alpha_n^{-1}\psi_n\circ\hat{w}_n\circ(2n\gamma_n)\rightarrow \psi\circ \hat{w}
\end{equation}
in $C([0,1],\ell^1)$ . The necessary tightness can be proved as follows:
\begin{eqnarray*}
\lefteqn{\lim_{\delta\rightarrow 0}\limsup_{n\rightarrow\infty}\sup_{\substack{s,t\in[0,1]\\|s-t|\leq \delta}}d_{\ell^{1}}\left(\alpha_n^{-1}\psi_n\left(\hat{w}_n\left(2n\gamma_n(s)\right)\right),\alpha_n^{-1}\psi_n\left(\hat{w}_n\left(2n\gamma_n(t)\right)\right)\right)}\\
&=&\lim_{\delta\rightarrow 0}\limsup_{n\rightarrow\infty}\sup_{\substack{s,t\in[0,1]\\|s-t|\leq \delta}}d_{\alpha_n^{-1}w_n}\left(\gamma_n(s),\gamma_n(t)\right)\hspace{80pt}\\
&=&\lim_{\delta\rightarrow 0}\sup_{\substack{s,t\in[0,1]\\|s-t|\leq \delta}}d_{w}\left(s,t\right)\\
&=&0,
\end{eqnarray*}
where $d_{\alpha_n^{-1}w_n}$ is a distance on $[0,1]$ defined similarly to the distance $d_w$ introduced at (\ref{dwdef}). The second equality is a result of Assumption \ref{a1}, and the final equality holds because $w$ is a continuous function. Now, by (\ref{mundef}), we can write
\[\tilde{\mu}_n=\lambda^{[0,1]}\circ (2n\gamma_n)^{-1}\circ \hat{w}_n^{-1}\circ \psi_n^{-1},\]
and so (\ref{funcconv}) implies that $\tilde{\mu}_n(\alpha_n\cdot)$ converges to $\tilde{\mu}=\lambda^{[0,1]}\circ \hat{w}^{-1}\circ\psi^{-1}$ weakly as probability measures on $\ell^1$. Thus, because $\tilde{\mu}^{(k)}=\tilde{\mu}\circ\pi_k^{-1}$ and $\tilde{\mu}_n^{(k)}=\tilde{\mu}_n\circ\pi_k^{-1}$ (these expressions follow from (\ref{proj}) and (\ref{projn}) respectively), we obtain that
\[\tilde{\mu}^{(k)}_n(\alpha_n\cdot)\rightarrow\tilde{\mu}^{(k)}\]
weakly as probability measures on $\ell^1$ for each $k\geq 1$.

The vector and measure convergence of the previous paragraph allows us to apply Proposition \ref{main3} to deduce that the law of \[\left(X^{\mathbf{x}_n}_{t\lambda^{\mathbf{x}_n}(\mathcal{T}^{\mathbf{x}_n})},A^{\mathbf{x}_n,\nu_n}_{t\lambda^{\mathbf{x}_n}(\mathcal{T}^{\mathbf{x}_n})}\right)_{t\geq 0}\]
under $\mathbf{P}^{\mathbf{x}_n}_0$ converges to the law of
\[\left(X^{\mathbf{x}}_{t\lambda^{\mathbf{x}}(\mathcal{T}^{\mathbf{x}})},A^{\mathbf{x},\nu}_{t\lambda^{\mathbf{x}}(\mathcal{T}^{\mathbf{x}})}\right)_{t\geq 0}\]
under $\mathbf{P}_0^{\mathbf{x}}$, weakly as probability measures on the space $C(\mathbb{R}_+,\ell^1)\times C(\mathbb{R}_+,\mathbb{R}_+)$, where $\nu_n:=\tilde{\mu}^{(k)}_n(\alpha_n\cdot)$ and $\nu:=\tilde{\mu}^{(k)}$. One can readily check that the distribution of the limit is equal to the distribution of $(\tilde{X}^{\mathcal{T}(k),\lambda^{(k)}},\hat{A}^{(k)})$ under $\mathbf{P}_\rho^{\mathcal{T}(k),\lambda^{(k)}}$. To complete the proof, we will use a coupling argument to show how the pairs $(\tilde{J}^{n,k},\hat{A}^{n,k})$ and $(X^{\mathbf{x}_n},A^{\mathbf{x}_n,\nu_n})$ can be related.

The defining properties of Brownian motion on a dendrite imply that, under $\mathbf{P}_0^{\mathbf{x}_n}$, the process $(X^{\mathbf{x}_n}_{h^{n}(m)})_{m\geq 0}$, where $h^{n}(0):=0$ and
\[h^{n}(m):=\inf\left\{t \geq h^{n}(m-1):\:d_{\ell^1}\left(X^{\mathbf{x}_n}_{t},X^{\mathbf{x}_n}_{h^{n}(m-1)}\right)=\alpha_n^{-1}\right\},\]
has the same law as $(\alpha_n^{-1}\tilde{J}^{n,k}_m)_{m\geq 0}$ under $\mathbf{P}^{T_n}_\rho$. In view of this fact, for the remainder of the proof we will abuse notation slightly by identifying $J^{n,k}$ with the process $(\psi_n^{-1}(\alpha_nX^{\mathbf{x}_n}_{h^{n}(m)}))_{m\geq 0}$, and $L^{n,k}$ with its local times. Furthermore, by considering the path segments of $X^{\mathbf{x}_n}$ between the hitting times $(h^n(m))_{m\geq0}$, it is possible to check that the rescaled increments $\alpha_n^{-2} (h^n(m)-h^n(m-1))$ are independently and identically distributed as the hitting time of $\{\pm 1\}$ by a standard Brownian motion in $\mathbb{R}$ started from zero, which is a random variable with mean 1 and finite fourth moments. After applying a standard martingale estimate (\cite{Kallenberg}, Proposition 7.16, for example), it readily follows that, for any $t_0,\varepsilon>0$,
\begin{equation}\label{hnkasm}
\lim_{n\rightarrow\infty}{\mathbf{P}}_0^{\mathbf{x}_n}\left(\sup_{m\leq t_0\alpha_n^2\Lambda_n^{(k)}}\left|h^{n}(m)-\frac{m}{\alpha_n^2}\right|>\varepsilon\right)=0,
\end{equation}
(cf. the proof of \cite{Croydoncbp}, Lemma 4.2). Consequently, since the sequence $(\mathbf{P}_0^{\mathbf{x}_n})_{n\geq 1}$ is tight and $\lambda^{\mathbf{x}_n}(\mathcal{T}^{\mathbf{x}_n})=\Lambda_n^{(k)}$, we are able to deduce that, for $t_0,\varepsilon>0$,
\begin{equation}\label{a}
\lim_{n\rightarrow\infty}{\mathbf{P}}_0^{\mathbf{x}_n}\left(\sup_{t\in[0,t_0]}
\left|\alpha_n^{-1}\tilde{J}^{n,k}_{t\alpha_n^2\Lambda^{(k)}_n}-X^{\mathbf{x}_n}_{t\lambda^{\mathbf{x}_n}(\mathcal{T}^{\mathbf{x}_n})}\right|>\varepsilon\right)=0.
\end{equation}

For the related additive functionals, to prove that
\begin{equation}\label{b}
\lim_{n\rightarrow\infty}{\mathbf{P}}_0^{\mathbf{x}_n}\left(\sup_{t\in[0,t_0]}\left|(n\alpha_n)^{-1}\hat{A}^{n,k}_{t\alpha_n^{2}\Lambda_n^{(k)}} -A^{\mathbf{x}_n,\nu_n}_{t\lambda^{\mathbf{x}_n}(\mathcal{T}^{\mathbf{x}_n})}\right|>\varepsilon\right)=0,
\end{equation}
by applying the tightness results of Lemmas \ref{lntight} and \ref{lnktight}, it will suffice to prove that
\begin{equation}\label{localclose}
\lim_{n\rightarrow\infty}\sup_{\sigma\in {T}_n(k)}
{\mathbf{P}}_0^{\mathbf{x}_n}\left( \sup_{t\in[0,t_0]}\left|\alpha_n^{-1}L^{n,k}_{t\alpha_n^{2}\Lambda_n^{(k)}}(\sigma)-
L^{\mathbf{x}_n}_{t\lambda^{\mathbf{x}_n}(\mathcal{T}^{\mathbf{x}_n})}(\alpha_n^{-1}\psi_n(\sigma))\right|>\varepsilon\right)=0,
\end{equation}
where the definition of $L^{n,k}$ is extended to continuous time by linear interpolation. To demonstrate that this is the case requires a simple adaptation of \cite{Croydoncbp}, Lemma 4.7. In particular, for $\sigma\in T_n(k)$ denote by $(\varsigma_i)_{i\geq 1}$ the hitting times of $\sigma$ by ${J}^{n,k}$, and define $\eta_i:={L}^{\mathbf{x}_n}_{h^{n}(\varsigma_i+1)}(\alpha_n^{-1}\psi_n(\sigma))-
{L}^{\mathbf{x}_n}_{h^{n}(\varsigma_i)}(\alpha_n^{-1}\psi_n(\sigma))$, so that the $(\eta_i)_{i=1}^\infty$ is a sequence of independent random variables, each distributed as ${2Z}/{\alpha_n\mathrm{deg}_{n,k}(\sigma)}$, where $Z$ represents the local time at zero of a standard Brownian motion in $\mathbb{R}$ started from zero, evaluated at the hitting time of $\{\pm 1\}$. By conditioning on $L^{n,k}_{t_0\alpha_n^2\Lambda_n^{(k)}}(\sigma)$, applying the fact that $Z$ is a random variable with mean 1 and finite fourth moments and recalling (\ref{polybound}), it is possible to check that
\begin{equation}\label{ee1}
\lim_{n\rightarrow \infty}\sup_{\sigma\in {T}_n(k)}{\mathbf{P}}_0^{\mathbf{x}_n}\left( \sup_{m\leq t_0\alpha_n^2\Lambda_n^{(k)}}\left|\eta_1+\dots+\eta_{\mathrm{deg}_{n,k}(\sigma)L^{n,k}_{m}(\sigma)/2}-
\alpha_n^{-1}L^{n,k}_{m}(\sigma)\right|>\varepsilon\right)=0,
\end{equation}
(cf. \cite{Croydoncbp}, equation (40)). Now, if $J^{n,k}_m=\sigma$, then  $\eta_1+\dots+\eta_{\mathrm{deg}_{n,k}(\sigma)L^{n,k}_{m}(\sigma)/2}$ is equal to $L^{\mathbf{x}_n}_{h^{n}(m+1)}(\alpha_n^{-1}\psi_n(\sigma))$, otherwise the sum is equal to $L^{\mathbf{x}_n}_{h^{n}(m)}(\alpha_n^{-1}\psi_n(\sigma))$. From this and the fact that, conditional on $J^{n,k}_m=\sigma$, ${L}^{\mathbf{x}_n}_{h^{n}(m+1)}(\alpha_n^{-1}\psi_n(\sigma))-
{L}^{\mathbf{x}_n}_{h^{n}(m)}(\alpha_n^{-1}\psi_n(\sigma))$ is distributed as
${2Z}/{\alpha_n\mathrm{deg}_{n,k}(\sigma)}$, it is readily deduced that (\ref{ee1}) also holds when $\alpha_n^{-1}L^{n,k}_{m}(\sigma)$ is replaced by the continuous local time $L^{\mathbf{x}_n}_{h^{n}(m)}(\alpha_n^{-1}\psi_n(\sigma))$. These observations, together with the hitting time estimate of (\ref{hnkasm}) and the local time tightness results of Lemmas \ref{lntight} and \ref{lnktight}, imply (\ref{localclose}), which thereby establishes (\ref{b}). Finally, combining (\ref{a}) and (\ref{b}) with the convergence result of the first part of the proof yields the lemma.
\end{proof}

{\lem\label{atight} If Assumption \ref{a1} holds, then, for $t_0,\varepsilon>0$,
\[\lim_{k\rightarrow\infty}\limsup_{n\rightarrow\infty}\mathbf{P}^{T_n}_\rho\left(\left( n\alpha_n\right)^{-1}\sup_{m\leq t_0 \alpha_n^2 \Lambda_n^{(k)}}\left|A_m^{n,k}-\hat{A}^{n,k}_m\right|>\varepsilon\right)=0.\]}
\begin{proof} The proof of this result is an adaptation of that used to obtain \cite{Croydoncbp}, Proposition 5.2. In particular, for $m\in\mathbb{N}$, we have that $|{A}^{n,k}_{m}-\hat{A}^{n,k}_{m}|$ is bounded above by
\[\left|\sum_{l=0}^{m-1}\left(A_{l+1}^{n,k}-A_l^{n,k}-\frac{2n\mu_n^{(k)}(\{J_l^{n,k}\})-2+\mathrm{deg}_{n,k}(J_l^{n,k})}{\mathrm{deg}_{n,k}(J_{l}^{n,k})}\right)\right|
+\sum_{l=0}^{m-1}\frac{|2-\mathrm{deg}_{n,k}(J_l^{n,k})|}{\mathrm{deg}_{n,k}(J_l^{n,k})}.\]
Clearly the second term is no greater than $m$, and so its supremum over $m\leq t_0 \alpha_n^2 \Lambda_n^{(k)}$ converges to 0 when multiplied by $(n\alpha_n)^{-1}$ (recall that $\alpha_n=o(n)$ by assumption). We now deal with the first term, which we will denote $\Sigma(m)$. Since, conditional on knowing $J^{n,k}$, the expected value of $A_{l+1}^{n,k}-A_l^{n,k}$ is precisely
\[\frac{2n\mu_n^{(k)}(\{J_l^{n,k}\})-2+\mathrm{deg}_{n,k}(J_l^{n,k})}{\mathrm{deg}_{n,k}(J_{l}^{n,k})}\]
and its expected square is bounded by
\[36n^2(\deg_{n,k}(J_l^{n,k})+\Delta_n^{(k)})\frac{\mu_n^{(k)}(\{J_l^{n,k}\})^2}{\mathrm{deg}_{n,k}(J_l^{n,k})},\]
where $\Delta_n^{(k)}:=\sup_{\sigma\in T_n}d_{T^n}(\sigma,\phi_{T_n,T_n(k)}(\sigma))$ (these are elementary simple random walk estimates, see \cite{Croydoncbp}, Lemma B.3), we can use Kolmogorov's maximum inequality (see \cite{Kallenberg}, Lemma 4.15) to deduce that, for $\varepsilon>0$,
\begin{eqnarray*}
\lefteqn{\mathbf{P}^{T_n}_\rho\left((n\alpha_n)^{-1}\sup_{m\leq t_0\alpha_n^2\Lambda_n^{(k)}}\Sigma(m)>\varepsilon\:\vline\: J^{n,k}\right)}\\
 &\leq &\frac{1}{n^{2}\alpha_n^2\varepsilon^2}\sum_{l=0}^{\lfloor t_0\alpha_n^2\Lambda_n^{(k)}\rfloor-1}36n^2(\deg_{n,k}(J_l^{n,k})+\Delta_n^{(k)})\frac{\mu_n^{(k)}(\{J_l^{n,k}\})^2}{\mathrm{deg}_{n,k}(J_l^{n,k})},\\
 &\leq&\frac{18}{n\alpha_n^2\varepsilon^2}\hat{A}^{n,k}_{\lfloor t_0\alpha_n^2\Lambda_n^{(k)}\rfloor}\left(\max_{\sigma\in T_n(k)}\mathrm{deg}_{n,k}(\sigma)+\Delta_n^{(k)}\right).
 \end{eqnarray*}
Hence, if we can show that
\begin{equation}\label{zero1}
\lim_{k\rightarrow\infty}\limsup_{n\rightarrow\infty}\alpha_n^{-1}\left(\max_{\sigma\in T_n(k)}\mathrm{deg}_{n,k}(\sigma)+\Delta_n^{(k)}\right)=0,
\end{equation}
and also
\begin{equation}\label{zero2}
\lim_{t\rightarrow\infty}\limsup_{k\rightarrow\infty}\limsup_{n\rightarrow\infty}
\mathbf{P}^{T_n}_\rho\left(\left(n\alpha_n\right)^{-1}\hat{A}^{n,k}_{t_0\alpha_n^{2}\Lambda_n^{(k)}}>t\right)=0,
\end{equation}
then the lemma will follow. That \begin{equation}\label{deltankdecay}
\lim_{k\rightarrow\infty}\limsup_{n\rightarrow\infty}\alpha_n^{-1}\Delta_n^{(k)}=0
\end{equation}
is a straightforward consequence of Assumption \ref{a1} (cf. \cite{Croydoncbp}, Lemma 2.7). Moreover, we clearly have that $\max_{\sigma\in T_n(k)}\mathrm{deg}_{n,k}(x)\leq k+1$, and so (\ref{zero1}) holds. The distributional convergence result of Lemma \ref{jumpconv} and the tail bound of (\ref{diverge}) together imply (\ref{zero2}), which completes the proof.
\end{proof}

We can now prove the convergence of simple random walks. In the statement of the result and the proof, discrete time processes are extended to continuous time by linear interpolation in the appropriate spaces.

{\prop\label{pconv} Under Assumption \ref{a1}, the laws of
\[\left(\alpha_n^{-1}\tilde{X}^{n}_{tn\alpha_n}\right)_{t\geq 0}\]
under $\mathbf{P}_\rho^{T_n}$ converge to $\tilde{\mathbf{P}}_\rho^{\mathcal{T},\mu}$ weakly as probability measures on the space $C(\mathbb{R}_+,\ell^1)$.}
\begin{proof} By Lemma \ref{akprops}, $\hat{A}^{(k)}$ is $\mathbf{P}_\rho^{\mathcal{T}(k),\lambda^{(k)}}$-a.s. continuous and strictly increasing. Consequently Proposition \ref{jumpconv} implies that the joint laws of the pairs
\[\left(\alpha_n^{-1}\tilde{J}^{n,k}_{t\alpha_n^2\Lambda_n^{(k)}}, (\alpha_n^{2}\Lambda_n^{(k)})^{-1}\hat{\tau}^{n,k}_{tn\alpha_n}\right)_{t\geq 0}\]
under $\mathbf{P}_\rho^{T_n}$ converge to the law of
\[\left(\tilde{X}^{\mathcal{T}(k),\lambda^{(k)}}_t,\hat{\tau}^{(k)}_t\right)_{t\geq 0}\]
under $\mathbf{P}_\rho^{\mathcal{T}(k),\lambda^{(k)}}$ weakly as probability measures on the space $C(\mathbb{R}_+,\ell^1)\times C(\mathbb{R}_+,\mathbb{R}_+)$. Hence we obtain from the continuous mapping theorem that the laws under $\mathbf{P}_\rho^{T_n}$ of $(\alpha_n^{-1}\psi_n(\hat{X}_{tn\alpha_n}^{n,k}))_{t\geq 0}$ converge to $\tilde{\mathbf{P}}_\rho^{\mathcal{T}(k),\mu^{(k)}}$ weakly as probability measures on the space $C(\mathbb{R}_+,\ell^1)$ for each $k\geq 1$, where we apply the representations of $\hat{X}^{n,k}$ and $X^{\mathcal{T}(k),\mu^{(k)}}$ from (\ref{hatx}) and Lemma \ref{timechange} respectively. Moreover, it is immediate from Proposition \ref{prockconv} and the continuity of $\psi$ that
$\tilde{\mathbf{P}}_\rho^{\mathcal{T}(k),\mu^{(k)}}$ converges to $\tilde{\mathbf{P}}_\rho^{\mathcal{T},\mu}$ as $k\rightarrow\infty$. The proposition will follow from these convergence results by applying \cite{Bill2}, Theorem 3.2, for example, if we can demonstrate the following tightness condition: for every $t_0,\varepsilon>0$,
\[\lim_{k\rightarrow\infty}\limsup_{n\rightarrow\infty}\mathbf{P}_\rho^{T_n}\left(\sup_{m\leq t_0n\alpha_n}
\alpha_n^{-1}d_{T_n}\left(X^n_{m},\hat{X}^{n,k}_{m}\right)>\varepsilon\right)=0.\]
To prove this, first observe that by construction $d_{T_n}(X^n_{m},{X}^{n,k}_{m})\leq \Delta_n^{(k)}$ for every $m\geq 0$, where $\Delta_n^{(k)}$ was defined in the proof of Lemma \ref{atight}. Thus (\ref{deltankdecay}) implies that it will actually suffice to prove that
\[\lim_{k\rightarrow\infty}\limsup_{n\rightarrow\infty}\mathbf{P}_\rho^{T_n}\left(\sup_{m\leq t_0n\alpha_n}
\alpha_n^{-1}d_{T_n}\left(X^{n,k}_{m},\hat{X}^{n,k}_{m}\right)>\varepsilon\right)=0.\]
Applying the representations of $X^{n,k}$ and $\hat{X}^{n,k}$ in terms of the jump-chain from (\ref{xnk}) and (\ref{hatx}) respectively, for any $\delta>0$ the probability in the left-hand side of this expression is bounded above by $p_1(n,k)+p_2(n,k)$, where
\[p_1(n,k):={\mathbf{P}_\rho^{T_n}\left(\sup_{\substack{s,t\in[0,t_0+\delta]\\|s-t|\leq \delta}}
\alpha_n^{-1}d_{\ell^1}\left(\tilde{J}^{n,k}_{\hat{\tau}^{n,k}(sn\alpha_n)},\tilde{J}^{n,k}_{\hat{\tau}^{n,k}(tn\alpha_n)}\right)>\varepsilon\right)}\]
\begin{eqnarray*}
\lefteqn{p_2(n,k):=}\\
&&\mathbf{P}_\rho^{T_n}\left(
{\tau}^{n,k}(tn\alpha_n)\not\in[\hat{\tau}^{n,k}((t-\delta)n\alpha_n\vee 0),\hat{\tau}^{n,k}((t+\delta)n\alpha_n)]\mbox{ for some }t\in[0,t_0]\right).
\end{eqnarray*}
It is elementary to check from Lemma \ref{atight} that
\[\lim_{k\rightarrow\infty}\limsup_{n\rightarrow\infty}p_2(n,k)=0\]
for any $\delta>0$. Furthermore, the convergence results of above yield
\[\lim_{k\rightarrow\infty}\limsup_{n\rightarrow\infty}p_1(n,k)=\mathbf{P}_\rho^{\mathcal{T},\mu}\left(\sup_{\substack{s,t\in[0,t_0+\delta]\\|s-t|\leq \delta}}d_{\mathcal{T}}\left(X^{\mathcal{T},\mu}_s,X^{\mathcal{T},\mu}_t\right)>\varepsilon\right),\]
which can be made arbitrarily small by letting $\delta\rightarrow 0$, because $X^{\mathcal{T},\mu}$ is a diffusion under $\mathbf{P}^{\mathcal{T},\mu}_\rho$. This completes the proof.
\end{proof}

We can now complete the proof of our main results.

\begin{proof}[Proof of Theorem \ref{qthm}]
Let $w\in\mathcal{W}$ satisfy $(\mathcal{T},\mu)\in\mathbb{T}^*$, and $U=(U_i)_{i\geq 1}$ be an independent sequence of $U[0,1]$ random variables. Since the Lebesgue measure $\lambda^{[0,1]}$ on $[0,1]$ has full support, and the measure $\mu$ is non-atomic, has full support and is supported on the leaves of $\mathcal{T}$, it is clear that $(U_i)_{i\geq1}$ is dense in $[0,1]$, and the vertices $(\sigma_i)_{i\geq 1}=(\hat{w}(U_i))_{i\geq1}$ are a dense collection of leaves of $\mathcal{T}$, distinct and not equal to $\rho$ for any $i$, for almost-every realisation of $U$. In particular, there exists a $u\in[0,1]^\mathbb{N}$ such that $(w,u)\in\Gamma$. Set $u^n=u$ for each $n$, so that under the assumptions of the theorem, we have that
\[\left(\alpha_n^{-1}w_n,u^n\right)\rightarrow(w,u)\]
in $C([0,1],\mathbb{R}_+)\times C(\mathbb{R}_+,\mathbb{R}_+)$ for some $(w,u)\in\Gamma$, which is Assumption \ref{a1}. Applying Proposition \ref{pconv}, this implies that
\[\tilde{\mathbf{P}}^{{T}_n}_\rho(\{f\in C([0,1],\ell^1):\:(\alpha_n^{-1}f(tn\alpha_n))_{t\in[0,1]}\in\cdot\})\rightarrow\tilde{\mathbf{P}}_\rho^{\mathcal{T},\mu}\]
weakly as probability measures on $C([0,1],\ell^1)$. The convergence $\tilde{\mu}_n(\alpha_n\cdot)\rightarrow\tilde{\mu}$ in $\mathcal{M}(\ell^1)$ was established in the proof of Lemma \ref{jumpconv}. That $\alpha_n^{-1}\tilde{T}_n(k)$ converges to $\tilde{\mathcal{T}}(k)$ in $\mathcal{K}(\ell^1)$ is a straightforward consequence of the convergence $\mathbf{x}_n\rightarrow\mathbf{x}$ in $(\ell^1)^k$, where the vertices $\mathbf{x}_n$ and $\mathbf{x}$ are defined as in the proof of Lemma \ref{jumpconv}. To extend this to the result that $\alpha_n^{-1}\tilde{T}_n\rightarrow\tilde{\mathcal{T}}$, we apply the tightness result of (\ref{deltankdecay}) and the fact that $\sup_{\sigma\in\mathcal{T}}d_\mathcal{T}(\sigma,\phi_{\mathcal{T},\mathcal{T}(k)}(\sigma))\rightarrow 0$, which was noted below (\ref{projproj}).
\end{proof}

\begin{proof}[Proof of Theorem \ref{athm}]
We start our proof, which is an adaptation of \cite{Croydoncbp}, Section 8, by outlining the construction of the measure $\mathbb{P}$. First, by following the proof of Lemma \ref{jumpconv}, it is possible to check that for each $k\geq 1$ the map $(w,u)\mapsto(\mathbf{x},\nu)$, where $\mathbf{x}:=(\psi(\sigma_1),\dots,\psi(\sigma_k))$ and $\nu:=\tilde{\mu}^{(k)}$, is continuous on $\Gamma$. By Lemma \ref{timechange}, Lemma \ref{akprops} and Proposition \ref{main3}, this implies that $(w,u)\mapsto (\tilde{T}(k),\tilde{\mu}^{(k)},\tilde{\mathbf{P}}^{\mathcal{T}(k),\mu^{(k)}}_\rho)$ is also continuous on $\Gamma$. Consequently, we obtain from Proposition \ref{prockconv} that  $(w,u)\mapsto (\tilde{T},\tilde{\mu},\tilde{\mathbf{P}}^{\mathcal{T},\mu}_\rho)$ is a measurable map on $\Gamma$. Given this, and noting that under the assumptions of the theorem that $(W,U)\in\Gamma$, $\mathbf{P}$-a.s., checking the existence of a unique probability measure satisfying (\ref{pdef}) is straightforward. The construction of $\mathbb{P}_n$ is similar, but easier. Finally, to check that $\mathbb{P}_n\circ\Theta_n^{-1}\rightarrow\mathbb{P}$, we consider a Skorohod-type coupling of random variables. In particular, since the relevant spaces are separable, if $(\alpha_n^{-1}W_n)_{n\geq 1}$ converges in distribution to $W$, then there exists a probability space upon which random variables $(W_n^*,U_n^*)_{n\geq 1}$ and $(W^*,U^*)$ are defined in such a way that $(W_n^*,U_n^*)$ has the distribution of $(W_n,U)$ for each $n$, $(W^*,U^*)$ has the distribution of $(W,U)$ and $(\alpha_n^{-1}W_n^*,U_n^*)\rightarrow (W^*,U^*)$, almost-surely. Defining all the related objects on this probability space, then exactly as in the previous proof we are able to deduce from Proposition \ref{pconv} that $\Theta_n(\tilde{T}_n,\tilde{\mu}_n,\tilde{\mathbf{P}}^{T_n}_\rho)\rightarrow (\tilde{T},\tilde{\mu},\tilde{\mathbf{P}}^{\mathcal{T},\mu}_\rho)$, almost-surely, and the result easily follows.
\end{proof}

\section{Application to $\alpha$-stable trees}\label{stable}

In this section, we describe the application of our results to Galton-Watson trees with a possibly infinite offspring distribution variance. Suppose $\xi$ is a non-negative integer-valued random variable that is aperiodic, has mean one, and is in the domain of attraction of a stable law with index $\alpha\in(1,2)$, by which we mean that there exists a sequence $a_n\uparrow \infty$ such that
\[\frac{\xi[n]-n}{a_n}\buildrel{d}\over{\rightarrow} \Xi,\]
where $\xi[n]$ is the sum of $n$ independent copies of $\xi$ and the limit random variable satisfies $\mathbf{E}(e^{-\lambda \Xi})=e^{-\lambda^\alpha}$. If $(T_n)_{n\geq 1}$ is a family of random trees such that $T_n$ is a Galton-Watson tree with offspring distribution $\xi$, conditioned on the total progeny being equal to $n$, then it is known (see \cite{Duqap}, Theorem 3.1) that the associated rescaled search-depth processes $(n^{-1}a_nW_n)_{n\geq 1}$ converge in distribution to a random excursion $W$ with law $N_{\alpha}^{(1)}$, say. The corresponding random real tree $\mathcal{T}$ is the $\alpha$-stable tree conditioned to have total mass equal to one, and we denote its law by $\Theta_\alpha^{(1)}$ (see also \cite{rrt}, Section 4). Consequently, to allow us to apply Theorem \ref{athm} to the sequence $(T_n)_{n\geq 1}$ with $\alpha_n=na_n^{-1}$, it remains to check that $\Theta_{\alpha}^{(1)}$-a.e. realisation of $\mu$ is non-atomic, supported on the leaves of $\mathcal{T}$ and satisfies (\ref{mulower}) for some $\kappa>0$. In fact, rather than investigating $\mu$ under the conditioned measure $\Theta_{\alpha}^{(1)}$, by rescaling it will suffice to check that the required properties hold under the unconditioned `excursion' measure $\Theta_\alpha$ of the $\alpha$-stable tree (see \cite{rrt}, Definition 4.2), as is done in the following proposition. For a rooted real tree, we use the notation $\chi(\mathcal{T}):=\sup\{d_\mathcal{T}(\rho,\sigma):\sigma\in\mathcal{T}\}$ to represent its height.

{\prop Let $\alpha\in(1,2)$. For $\Theta_\alpha$-a.e. realisation of $\mathcal{T}$, $\mu$ is non-atomic, supported on the leaves of $\mathcal{T}$ and satisfies
\begin{equation}\label{loglower}
\liminf_{r\rightarrow 0}\frac{\inf_{\sigma\in\mathcal{T}}\mu(B(\sigma,r))}{r^{\frac{\alpha}{\alpha-1}}\left(\ln r^{-1}\right)^{-\frac{\alpha}{\alpha-1}}}>0.
\end{equation}
In particular, (\ref{loglower}) implies (\ref{mulower}) for any $\kappa>\alpha/(\alpha-1)$.}
\begin{proof} That $\mu$ is non-atomic and supported on the leaves of $\mathcal{T}$ for $\Theta_\alpha$-a.e. realisation of $\mathcal{T}$ follows from results of \cite{LegallDuquesne}. Thus it remains to check (\ref{loglower}).

Fix a compact rooted real tree $\mathcal{T}$ and $r>0$. Following \cite{LegallDuquesne}, denote by $(\mathcal{T}^{(i),o})_{i\in\mathcal{I}}$ the connected components of the open set $\{\sigma\in\mathcal{T}:d_\mathcal{T}(\rho,\sigma)>r\}$. Define $\mathcal{T}^{(i)}:=\mathcal{T}^{(i),o}\cup\{\sigma_i\}$, where $\sigma_i$ is the common ancestor in level $r$ of every $\sigma\in\mathcal{T}^{(i),o}$, so that $\mathcal{T}^{(i)}$ is a compact rooted real tree, and we set its root to be $\sigma_i$. The trees $(\mathcal{T}^{(i)})_{i\in\mathcal{I}}$ are the subtrees of $\mathcal{T}$ originating from level $r$. If $\chi(\mathcal{T}^{(i)})\geq\delta$, then we say that $\mathcal{T}^{(i)}$ hits level $r+\delta$. The number of subtrees of $\mathcal{T}$ originating from level $r$ that hit level $r+\delta$ will be written $Z(r,\delta)$.

For integers $n,k\geq 0$, define $(\mathcal{T}^{n,k,(i)})_{i=1}^{Z(k2^{-n},2^{-n})}$ to be the collection of subtrees of $\mathcal{T}$ originating at level $k2^{-n}$ that hit level $(k+1)2^{-n}$. If we assume that $\chi(\mathcal{T})\geq 2^{-n}$, then it is elementary to check that, for fixed $r>0$,
\[\inf_{\sigma\in B(\rho,r)}\mu\left(B(\sigma,3.2^{-n})\right)\geq \inf_{\substack{0\leq k\leq 2^n r\\ 1\leq i\leq Z(k2^{-n},2^{-n})}}\mu \left(\mathcal{T}^{n,k,(i)}(2^{-n})\right),\]
where $\mathcal{T}^{n,k,(i)}(2^{-n})$ is the ball in $\mathcal{T}^{n,k,(i)}$ of radius $2^{-n}$ centred at the root of $\mathcal{T}^{n,k,(i)}$, which implies that
\begin{eqnarray}
\lefteqn{\Theta_\alpha\left(\inf_{\sigma\in B(\rho,r)}\mu\left(B(\sigma,3.2^{-n})\right)\leq x,\chi(\mathcal{T})\geq 2^{-n}\right)}\nonumber\\
&\leq& \sum_{0\leq k\leq 2^n r}\Theta_\alpha\left(\inf_{1\leq i\leq Z(k2^{-n},2^{-n})}\mu \left(\mathcal{T}^{n,k,(i)}(2^{-n})\right)\leq x,\chi(\mathcal{T})\geq (k+1)2^{-n}\right).\label{sum}
\end{eqnarray}
Note that if $\chi(\mathcal{T})< (k+1)2^{-n}$ for some $k\geq 0$, then $Z(k2^{-n},2^{-n})=0$ and the infimum appearing in the $k$th summand is infinite. Hence the summands are not decreased by including the statement $\chi(\mathcal{T})\geq (k+1)2^{-n}$ as we do.

The branching property of L\'{e}vy trees (\cite{LegallDuquesne}, Theorem 4.2) implies that under the measure $\Theta_\alpha(\cdot|\chi(\mathcal{T})\geq k2^{-n})$, conditional on $Z(k2^{-n},2^{-n})$, the trees in the collection $(\mathcal{T}^{n,k,(i)})_{i=1}^{Z(k2^{-n},2^{-n})}$ are distributed as an independent sample chosen according to the law $\Theta_\alpha(\cdot|\chi(\mathcal{T})\geq 2^{-n})$. Hence, writing $\tilde{\Theta}_\alpha:=\Theta_\alpha(\cdot|\chi(\mathcal{T})\geq k2^{-n})$ and $Z:=Z(k2^{-n},2^{-n})$, when $k\geq 1$ the $k$th summand of (\ref{sum}) satisfies
\begin{eqnarray*}
\lefteqn{\Theta_\alpha\left(\inf_{1\leq i\leq Z}\mu \left(\mathcal{T}^{n,k,(i)}(2^{-n})\right)\leq x,\chi(\mathcal{T})\geq (k+1)2^{-n}\right)}\\
&=&\tilde{\Theta}_\alpha\left(\tilde{\Theta}_\alpha\left(\inf_{1\leq i\leq Z}\mu \left(\mathcal{T}^{n,k,(i)}(2^{-n})\right)\leq x,\chi(\mathcal{T})\geq (k+1)2^{-n}|Z\right)\right)\\
&&\times\Theta_\alpha\left(\chi(\mathcal{T})\geq k2^{-n}\right)\\
&=&\tilde{\Theta}_\alpha\left(\mathbf{1}_{\{Z\neq 0\}}\tilde{\Theta}_\alpha\left(\inf_{1\leq i\leq Z}\mu \left(\mathcal{T}^{n,k,(i)}(2^{-n})\right)\leq x|Z\right)\right)\Theta_\alpha\left(\chi(\mathcal{T})\geq k2^{-n}\right)\\
&\leq & \tilde{\Theta}_\alpha\left(Z{\Theta_\alpha}\left(\mu (B(\rho,2^{-n}))\leq x|\chi(\mathcal{T})\geq 2^{-n}\right)\right)\Theta_\alpha\left(\chi(\mathcal{T})\geq k2^{-n}\right)\\
&=&{\Theta_\alpha}\left(Z\right) {\Theta_\alpha}\left(\mu( B(\rho,2^{-n}))\leq x|\chi(\mathcal{T})\geq 2^{-n}\right).
\end{eqnarray*}
Moreover, we have that ${\Theta_\alpha}(Z)={\Theta_\alpha}(\chi(\mathcal{T})\geq 2^{-n})$ (cf. proof of \cite{LegallDuquesne}, Lemma 5.4), and so summing over $k$ yields
\begin{eqnarray*}
\lefteqn{\Theta_\alpha\left(\inf_{\sigma\in B(\rho,r)}\mu\left(B(\sigma,3.2^{-n})\right)\leq x,\chi(\mathcal{T})\geq 2^{-n}\right)}\\
&\leq& (2^nr+1){\Theta_\alpha}\left(\mu( B(\rho,2^{-n}))\leq x,\chi(\mathcal{T})\geq 2^{-n}\right)\\
&\leq& c_1 (2^nr+1) 2^{\frac{n}{\alpha-1}}e^{-c_22^{-n}x^{-\frac{\alpha-1}{\alpha}}},
\end{eqnarray*}
where $c_1$ and $c_2$ are constants not depending on $n$, $r$ or $x$, and the second inequality is a result of equation (5.8) of \cite{treemeas}. Taking $c_3$ suitably small, this implies for any $R>0$ that
\[\sum_{n=0}^{\infty}\Theta_\alpha\left(\inf_{\sigma\in B(\rho,r)}\mu\left(B(\sigma,3.2^{-n})\right)\leq c_3(n2^n)^{-\frac{\alpha}{\alpha-1}},\chi(\mathcal{T})\geq R\right)<\infty.\]
The result follows by applying the Borel-Cantelli lemma, and then letting $r\rightarrow \infty$ and $R\rightarrow 0$.
\end{proof}

To complete this section, we note that the above proposition allows the heat kernel estimate of Lemma \ref{est}(b) to be improved in the following way in the $\alpha$-stable case. By comparison with results appearing in \cite{CK} for random walks on infinite variance Galton-Watson trees conditioned to survive, we expect that the polynomial term in the following lemma is the best possible. For an analogous estimate in the Brownian case, see \cite{Croydoncrt}.

{\cor  Let $\alpha\in(1,2)$. For $\Theta_\alpha$-a.e. realisation of $\mathcal{T}$, the Brownian motion $X^{\mathcal{T},\mu}$ admits a transition density $(p_t(\sigma,\sigma'))_{\sigma,\sigma'\in\mathcal{T},t>0}$ that satisfies
\[\limsup_{t\rightarrow 0}\frac{\sup_{\sigma,\sigma'\in\mathcal{T}}p_t(\sigma,\sigma')}{t^{-\frac{\alpha}{2\alpha-1}}\left(\ln t^{-1}\right)^{\frac{\alpha}{2\alpha-1}}}<\infty.\]}
\def\cprime{$'$}
\providecommand{\bysame}{\leavevmode\hbox to3em{\hrulefill}\thinspace}
\providecommand{\MR}{\relax\ifhmode\unskip\space\fi MR }
\providecommand{\MRhref}[2]{%
  \href{http://www.ams.org/mathscinet-getitem?mr=#1}{#2}
}
\providecommand{\href}[2]{#2}


\begin{thebibliography}{10}

\bibitem{Aldous2}
D.~Aldous, \emph{The continuum random tree. {II}. {A}n overview}, Stochastic
  analysis (Durham, 1990), London Math. Soc. Lecture Note Ser., vol. 167,
  Cambridge Univ. Press, Cambridge, 1991, pp.~23--70.

\bibitem{Aldous3}
\bysame, \emph{The continuum random tree. {III}}, Ann. Probab. \textbf{21}
  (1993), no.~1, 248--289.

\bibitem{Bill2}
P.~Billingsley, \emph{Convergence of probability measures}, second ed., Wiley
  Series in Probability and Statistics: Probability and Statistics, John Wiley
  \& Sons Inc., New York, 1999, A Wiley-Interscience Publication.

\bibitem{BlumGet}
R.~M. Blumenthal and R.~K. Getoor, \emph{Markov processes and potential
  theory}, Pure and Applied Mathematics, Vol. 29, Academic Press, New York,
  1968.

\bibitem{Borodin}
A.~N. Borodin, \emph{The asymptotic behavior of local times of recurrent random
  walks with finite variance}, Teor. Veroyatnost. i Primenen. \textbf{26}
  (1981), no.~4, 769--783.

\bibitem{Croydon}
D.~A. Croydon, \emph{Heat kernel fluctuations for a resistance form with
  non-uniform volume growth}, Proc. Lond. Math. Soc. (3) \textbf{94} (2007),
  no.~3, 672--694.

\bibitem{Croydoncbp}
\bysame, \emph{Convergence of simple random walks on random discrete trees to
  {B}rownian motion on the continuum random tree}, Ann. Inst. Henri Poincar\'e
  Probab. Stat. \textbf{44} (2008), no.~6, 987--1019.

\bibitem{Croydoncrt}
\bysame, \emph{Volume growth and heat kernel estimates for the continuum random
  tree}, Probab. Theory Related Fields \textbf{140} (2008), no.~1-2, 207--238.

\bibitem{CHLLT}
D.~A. Croydon and B.~M. Hambly, \emph{Local limit theorems for sequences of
  simple random walks on graphs}, Potential Anal. \textbf{29} (2008), no.~4,
  351--389.

\bibitem{CK}
D.~A. Croydon and T.~Kumagai, \emph{Random walks on {G}alton-{W}atson trees
  with infinite variance offspring distribution conditioned to survive},
  Electron. J. Probab. \textbf{13} (2008), no. 51, 1419--1441.

\bibitem{DudleyGauss}
R.~M. Dudley, \emph{Sample functions of the {G}aussian process}, Ann.
  Probability \textbf{1} (1973), no.~1, 66--103.

\bibitem{Duqap}
T.~Duquesne, \emph{A limit theorem for the contour process of conditioned
  {G}alton-{W}atson trees}, Ann. Probab. \textbf{31} (2003), no.~2, 996--1027.

\bibitem{LegallDuquesne}
T.~Duquesne and J.-F. Le~Gall, \emph{Probabilistic and fractal aspects of
  {L}\'evy trees}, Probab. Theory Related Fields \textbf{131} (2005), no.~4,
  553--603.

\bibitem{treemeas}
\bysame, \emph{The {H}ausdorff measure of stable trees}, ALEA \textbf{1}
  (2006), 393--415.

\bibitem{FOT}
M.~Fukushima, Y.~{\=O}shima, and M.~Takeda, \emph{Dirichlet forms and symmetric
  {M}arkov processes}, de Gruyter Studies in Mathematics, vol.~19, Walter de
  Gruyter \& Co., Berlin, 1994.

\bibitem{Kallenberg}
O.~Kallenberg, \emph{Foundations of modern probability}, second ed.,
  Probability and its Applications (New York), Springer-Verlag, New York, 2002.

\bibitem{Kigamidendrite}
J.~Kigami, \emph{Harmonic calculus on limits of networks and its application to
  dendrites}, J. Funct. Anal. \textbf{128} (1995), no.~1, 48--86.

\bibitem{Kumagai}
T.~Kumagai, \emph{Heat kernel estimates and parabolic {H}arnack inequalities on
  graphs and resistance forms}, Publ. Res. Inst. Math. Sci. \textbf{40} (2004),
  no.~3, 793--818.

\bibitem{rrt}
J.-F. Le~Gall, \emph{Random real trees}, Ann. Fac. Sci. Toulouse Math. (6)
  \textbf{15} (2006), no.~1, 35--62.

\bibitem{MarcusRosen}
M.~B. Marcus and J.~Rosen, \emph{Sample path properties of the local times of
  strongly symmetric {M}arkov processes via {G}aussian processes}, Ann. Probab.
  \textbf{20} (1992), no.~4, 1603--1684.

\end{thebibliography}
\end{document}